\documentclass[10pt,compsoc,journal]{IEEEtran}

\ifCLASSOPTIONcompsoc
  \usepackage[nocompress]{cite}
\else
  \usepackage{cite}
\fi
\usepackage{hyperref}
\usepackage{graphicx}
\usepackage{amsfonts}
\usepackage{amsmath}
\usepackage{amssymb}
\usepackage{amsthm}
\usepackage{float}
\usepackage{algorithm}
\usepackage[noend]{algpseudocode}
\usepackage{array}
\usepackage{listings}
\usepackage{subcaption}

\newtheorem{example}{Example}
\newtheorem{remark}{Remark}
\newtheorem{lemma}{Lemma}

\begin{document}

\title{Numerically Stable Recurrence Relations for the Communication Hiding Pipelined Conjugate Gradient Method}

\author{Siegfried~Cools$^*$,
        Jeffrey~Cornelis$^*$,
        Wim Vanroose$^*$
\IEEEcompsocitemizethanks{\IEEEcompsocthanksitem[$^*$] 
Applied Mathematics Group, Department of Mathematics and Computer Science, University of Antwerp.
\textbf{Address:} University of Antwerp, Campus Middelheim, Building G, Middelheimlaan 1, 2020 Antwerp, Belgium. 
\textbf{E-mail:} siegfried.cools@uantwerp.be (corresponding author). 
\textbf{Funding:} S.\,Cools is funded by Research Foundation Flanders (FWO) under grant 12H4617N. J.\,Cornelis receives funding from the University of Antwerp Research Council under the University Research Fund (BOF).}
\thanks{Manuscript received \today; revised (TBA).}}

\markboth{IEEE Transactions on Parallel and Distributed Systems,~Vol.~?, No.~?, February~2019}%
{Cools \MakeLowercase{\textit{et al.}}: Stable Recurrence Relations for Pipelined Conjugate Gradients}
%


\IEEEtitleabstractindextext{%
\begin{abstract}
Pipelined Krylov subspace methods (also referred to as communication-hiding methods) have been proposed in the literature as a scalable alternative to classic Krylov subspace algorithms for iteratively computing the solution to a large linear system in parallel. For symmetric and positive definite system matrices the pipelined Conjugate Gradient method, p($l$)-CG, outperforms its classic Conjugate Gradient counterpart on large scale distributed memory hardware by overlapping global communication with essential computations like the matrix-vector product, thus ``hiding'' global communication. A well-known drawback of the pipelining technique is the (possibly significant) loss of numerical stability. In this work a numerically stable variant of the pipelined Conjugate Gradient algorithm is presented that avoids the propagation of local rounding errors in the finite precision recurrence relations that construct the Krylov subspace basis. The multi-term recurrence relation for the basis vector is replaced by $\ell$ three-term recurrences, improving stability without increasing the overall computational cost of the algorithm. The proposed modification ensures that the pipelined Conjugate Gradient method is able to attain a highly accurate solution independently of the pipeline length. Numerical experiments demonstrate a combination of excellent parallel performance and improved maximal attainable accuracy for the new pipelined Conjugate Gradient algorithm. This work thus resolves one of the major practical restrictions for the useability of pipelined Krylov subspace methods.\vspace{-0.2cm}
\end{abstract}

\begin{IEEEkeywords}
Krylov subspace methods, Pipelining, Parallel performance, Global communication, Latency hiding, Conjugate Gradients, Numerical stability, Inexact computations, Attainable accuracy.\\ \vspace{-0.1cm} 
\end{IEEEkeywords}}
	
\maketitle
%

\IEEEraisesectionheading{\section{Introduction}\label{sec:introduction}}

\thispagestyle{plain}
\pagestyle{plain}

\IEEEPARstart{T}{he} 
family of iterative solvers known as Krylov subspace methods (KSMs) \cite{
liesen2012krylov,saad2003iterative,van2003iterative} are among the most efficient present-day methods for solving large scale sparse systems of linear equations. The mother of all Krylov subspace methods is undoubtedly the Conjugate Gradient method (CG) that was 
derived in 1952 \cite{hestenes1952methods} to the aim of solving linear systems $Ax=b$ with a symmetric positive definite and preferably sparse matrix $A$. The CG method is one of the most widely used methods for solving said systems today, which form the basis of a plethora of scientific and industrial applications. However, driven by the 
essential transition from optimal single node performance towards massively parallel computer hardware over the last decades \cite{shalf2007new}, the bottleneck for fast execution of Krylov subspace methods has shifted. Whereas in the past the application of the sparse matrix-vector product (\textsc{spmv}) was considered the most time-consuming part of the algorithm, the global synchronizations required in dot product and norm computations form the main bottleneck for efficient execution on present-day distributed memory hardware \cite{dongarra2015hpcg}.

Driven by the increasing levels of parallelism in present-day HPC machines, as attested by the current strive towards exascale high-performance computing software \cite{dongarra2011international}, research on the elimination of the global communication bottleneck has recently regained significant attention from the international computer science, engineering and numerical mathematics communities.
Sprouting largely from pioneering work on reducing communication in Krylov subspace methods from the late 1980's and 90's \cite{strakovs1987effectivity,meurant1987multitasking,chronopoulos1989s,d1992reducing,demmel1993parallel,erhel1995parallel
}, a number of variants of the classic Krylov subspace algorithms have been introduced over the last years. We point out recent work by Chronopoulos et al.~\cite{chronopoulos2010block}, Hoemmen \cite{hoemmen2010communication}, Carson et al.~\cite{carson2013avoiding}, McInnes et al.~\cite{mcinnes2014hierarchical}, Grigori et al.~\cite{grigori2016enlarged}, Eller et al.~\cite{eller2016scalable}, Imberti et al.~\cite{imberti2017varying} and Zhuang et al.~\cite{zhuang2017iteration}.

The contents of the 
 current work are situated in the research branch on so-called ``\emph{pipelined}'' Krylov subspace methods\footnote{\textbf{Note:} In the context of communication reduction in Krylov subspace methods, the terminology \emph{``pipelined''} KSMs that is used throughout the related applied linear algebra and computer science literature refers to \emph{software pipelining}, i.e. algorithmic reformulations to the KSM procedure in order to reduce communication overhead, and should not be confused with hardware-level \emph{instruction pipelining} (ILP).} \cite{ghysels2013hiding,ghysels2014hiding,cools2017communication}. Alternatively called ``\emph{communication hiding}'' methods, these algorithmic variations to classic Krylov subspace methods are designed to overlap the time-consuming global communications in each iteration of the algorithm with computational tasks such as calculating \textsc{spmv}s or \textsc{axpy}s (vector operations of the form $y \leftarrow \alpha x + y$). Thanks to the reduction/elimination of the synchronization bottleneck, pipelined algorithms have been shown to \emph{increase parallel efficiency} 
by allowing the algorithm to continue scaling on large numbers of processors \cite{sanan2016pipelined,yamazaki2017improving}. However, the algorithmic reformulations that allows for this efficiency increase come at the cost of \emph{reduced numerical stability} \cite{ghysels2014hiding,carson2018numerical}, which presently is one of the main drawbacks of pipelined (as well as other communication reducing) methods. Research on analyzing and improving the numerical stability of pipelined Krylov subspace methods, which is essential both for a proper understanding and the practical usability of the methods, has recently been performed by the authors \cite{cools2018analyzing,cools2018analyzing2} and others \cite{carson2014residual,carson2018numerical}.  

This work presents a numerically stable variant of the $l$-length pipelined Conjugate Gradient method, p($l$)-CG for short. The p($l$)-CG method was 
presented in \cite{cornelis2017communication} and allows to overlap each global reduction phase with the computational work of $l$ subsequent iterations. The pipeline length $l$ is a parameter of the method that can be chosen depending on the problem and hardware setup (as a function of the communication-to-computation ratio). As is the case for all communication reducing Krylov subspace methods, the preconditioner choice influences the communication-to-computation ratio and thus affects the performance of the method. The pipeline length hence also depends on the effort invested in the preconditioner. A preconditioner that uses limited global communication (block Jacobi, no-overlap DDM, \ldots) is generally preferred in this setting.

The propagation of local rounding errors in the multi-term recurrence relations of the p($l$)-CG algorithm is the primary source of loss of attainable accuracy on the final solution \cite{cornelis2017communication}. By introducing intermediate auxiliary basis variables, we derive a new p($l$)-CG algorithm with modified recurrence relations for which \emph{no rounding error propagation} occurs. It is proven analytically that the resulting recurrence relations are numerically stable for \emph{any} pipeline length $l$. 
The new 
algorithm is guaranteed to reach the same accuracy as the classic CG method. This work thus resolves one of the major restrictions for the practical use of pipelined Krylov subspace methods. The redesigned algorithm comes at no additional computational cost and has only a minor storage overhead compared to the former p($l$)-CG algorithm, thus effectively replacing the earlier implementation of the method. In addition, it is shown that formulating a preconditioned version of the new algorithm is straightforward. 

We conclude this introduction by providing a short overview of the further contents of this manuscript. Section \ref{sec:pipelcg} presents a high-level summary of the basic principles behind the $l$-length pipelined CG method and formulates the key numerical properties of the method that motivate this work. 
It 
familiarizes the reader with 
the notations and concepts 
used throughout this paper.
In Section \ref{sec:analyzing} the numerical stability analysis of the p($l$)-CG recurrence relations is briefly recapped, as it forms the basis for the analysis of the stable algorithm in Section \ref{sec:rounding}. Section \ref{sec:stabilizing} contains the main contributions of this work, presenting the technical details of the stable p($l$)-CG algorithm alongside an overview of its main implementation properties and a numerical analysis of the new rounding error resilient recurrence relations. Numerical experiments that demonstrate both the parallel performance of the p($l$)-CG method and the 
excellent attainable accuracy in comparison to earlier variants of pipelined Krylov subspace methods are presented in Section \ref{sec:numerical}. The manuscript is concluded in Section \ref{sec:conclusions}.

For completeness we note that the numerical analysis in Section \ref{sec:stabilizing} focuses on analyzing the propagation of local rounding errors throughout the new p($l$)-CG algorithm in detail, but does not include a standard forward or backward stability analysis with bounds on the 
local rounding errors.


\section{Deep pipelined Conjugate Gradients} \label{sec:pipelcg}

The deep pipelined Conjugate Gradient method, denoted p($l$)-CG for short, was first presented in \cite{cornelis2017communication}, where it was derived in analogy to the p($l$)-GMRES method \cite{ghysels2013hiding}. The parameter $l$ 
represents the pipeline length which indicates the number of iterations that are overlapped by each global reduction phase. We summarize the current state-of-the-art deep pipelined p($l$)-CG method 
below, 
which forms the starting point for the discussion in this work.

\subsection{Basis 
recurrence relations 
in exact arithmetic}

Let $V_{i-l+1}=[v_{0},v_{1},\ldots,v_{i-l}]$ be the orthonormal basis for the Krylov subspace $\mathcal{K}_{i-l+1}(A,v_{0})$ in iteration $i$ of the p($l$)-CG algorithm, consisting of $i-l+1$ vectors. Here $A$ is a symmetric positive definite matrix. The Krylov subspace basis vectors satisfy the Lanczos relation
\begin{equation} \label{eq:Arnoldi_v}
AV_{j} = V_{j+1} T_{j+1,j}, \qquad 1 \leq j \leq i-l.
\end{equation}
with \vspace{-0.3cm}
\begin{equation}
T_{j+1,j} =  
	\begin{pmatrix} 
  \gamma_{0} & \delta_{0} & & \\ 
 \delta_{0} & \gamma_{1} & \ddots & \\  
  & \ddots & \ddots & \delta_{j-2} \\ 
  & & \delta_{j-2} & \gamma_{j-1}  \\
  & & & \delta_{j-1}
	\end{pmatrix}.
\end{equation}
Let $\delta_{-1} = 0$, then the Lanczos relation \eqref{eq:Arnoldi_v} translates in vector notation to 
\begin{equation}\label{eq:v_ex}
  v_{j+1}= ( Av_{j} - \gamma_j v_j - \delta_{j-1} v_{j-1} ) / \delta_j, \quad 0 \leq j < i - l.
\end{equation}
The auxiliary basis $Z_{i+1}=[z_{0},z_{1},\ldots,z_{i}]$ runs $l$ vectors ahead of the basis $V_{i-l+1}$ and is defined as 
\begin{equation} \label{eq:z_ex}
  z_{j}= \left\{ \begin{matrix} v_{0}, & j=0, \\ P_{j}(A)v_{0}, & 0<j\leq l, \\ P_{l}(A)v_{j-l}, & l<j\leq i, \end{matrix} \right. 
\end{equation}
where the matrix polynomial $P_{l}(A)$ is given by
\begin{equation} \label{eq:poly}
	P_{k}(t) = \prod_{j=0}^{k-1} (t-\sigma_{j}), \qquad k\leq l,
\end{equation}
with optional 
stabilizing shifts $\sigma_{j}\in \mathbb{R}$, see \cite{ghysels2013hiding,cornelis2017communication,hoemmen2010communication}. 
We refer to Section \ref{sec:conditioning} for a discussion on the Krylov subspace basis $Z_{i+1}$, i.e.~the choice of the polynomial $P_l(A)$.
Contrary to the basis $V_{i-l+1}$, the auxiliary basis $Z_{i+1}$ is in general not orthonormal.
It is constructed using the recursive definitions \vspace{-0.5cm}
\begin{equation} \label{eq:z_rec}
  z_{j+1} = \left\{
	\begin{matrix}
			(A-\sigma_{j}I) \, z_{j}, & 0 \leq j < l, \\
			(Az_{j} - \gamma_{j-l} z_j - \delta_{j-l-1} z_{j-1} ) / \delta_{j-l}, & l \leq j < i, 
	\end{matrix} \right.
\end{equation}
which are obtained by multiplying the Lanczos relation \eqref{eq:v_ex} on both sides by $P_{l}(A)$.
Expression \eqref{eq:z_rec} translates into a Lanczos type matrix relation
\begin{equation}
A Z_j = Z_{j+1} B_{j+1,j}, \qquad 1 \leq j \leq i,
\end{equation}
where the matrix $B_{j+1,j}$ contains the matrix $T_{j-l+1,j-l}$, which is shifted $l$ places along the main diagonal.
The bases $V_{j}$ and $Z_j$ are connected through the basis transformation $Z_{j}=V_{j}G_{j}$ for $1 \leq j \leq i-l+1$, where $G_j$ is a banded upper triangular matrix with a band width of $2l+1$ non-zero diagonals \cite{cornelis2017communication}. For a symmetric matrix $A$ the matrix $G_{i+1}$ is symmetric around its $l$-th upper diagonal, since
\begin{align} \label{eq:symmetry_G}
	g_{j,i} &= (z_{i},v_{j}) = (P_{l}(A)v_{i-l},v_{j}) = (v_{i-l},P_{l}(A)v_{j}) \notag \\ 
					&= (v_{i-l},z_{j+l}) = g_{i-l,j+l}.
\end{align}
The following recurrence relation for $v_{j+1}$ is derived from the basis transformation (with $0 \leq j < i-l$):
\begin{equation} \label{eq:v_rec}
  v_{j+1} = \left( z_{j+1} - \sum_{k=j-2l+1}^{j} g_{k,j+1} v_{k} \right) / g_{j+1,j+1}.
\end{equation}
A total of $l$ iterations after the dot-products 
\begin{equation}
g_{j,i+1} = \left\{ \begin{matrix}
	(z_{i+1},v_{j}); & j=\max(0,i-2l+1),\ldots,i-l+1, \\ 
	(z_{i+1},z_{j});  &j=i-l+2,\ldots,i+1, \end{matrix}
	\right.
\end{equation}
have been initiated, the elements $g_{j,i-l+1}$ with $i-2l+2 \leq j \leq i-l +1$, which were computed as $(z_{i-l+1},z_{j})$, are corrected as follows:
\begin{equation}
g_{j,i-l+1} = \frac{g_{j,i-l+1}-\sum_{k=i-3l+1}^{j-1}g_{k,j}g_{k,i-l+1}}{g_{j,j}}, 
\end{equation}
for $j=i-2l+2,\ldots, i-l$, and:
\begin{equation}
g_{i-l+1,i-l+1} =\sqrt{g_{i-l+1,i-l+1}-\sum_{k=i-3l+1}^{i-l}g_{k,i-l+1}^2}.
\end{equation}
Additionally, in the $i$-th iteration the tridiagonal matrix $T_{i-l+2,i-l+1}$, see \eqref{eq:Arnoldi_v}, can be updated recursively by adding one column. 
The diagonal element $\gamma_{i-l}$ is characterized by the expressions: $\gamma_{i-l} = $
\begin{equation} \label{eq:haa1}
\left\{
  \begin{aligned} 
    & \frac{g_{i-l,i-l+1}+\sigma_{i-l}g_{i-l,i-l} -g_{i-l-1,i-l}\delta_{i-l-1}}{g_{i-l,i-l}}, \quad l \leq i < 2l,\\
    & \frac{g_{i-l,i-l}\gamma_{i-2l}+g_{i-l,i-l+1}\delta_{i-2l}-g_{i-l-1,i-l}\delta_{i-l-1}}{g_{i-l,i-l}}, ~i \geq 2l.
  \end{aligned} 
	\right.
\end{equation}
The term $-g_{i-l-1,i-l}\delta_{i-l-1}$
is considered zero when $i = l$.
The update for the off-diagonal element $\delta_{i-l}$ is given by
\begin{equation} \label{eq:haa3}
\delta_{i-l} = \left\{ 
  \begin{matrix} 
    g_{i-l+1,i-l+1}/g_{i-l,i-l}, & l \leq i < 2l, \\ 
    (g_{i-l+1,i-l+1}\delta_{i-2l})/g_{i-l,i-l}, & i \geq 2l.
  \end{matrix} 
	\right.
\end{equation}
The element $\delta_{i-l-1}$ has already been computed in the previous iteration and can thus simply be copied due to the symmetry of $T_{i-l+2,i-l+1}$. 

\begin{algorithm*}[t]
{\small
\caption{Original $l$-length pipelined p($l$)-CG \cite{cornelis2017communication} \hfill \textbf{Input:} $A$, $b$, $x_0$, $l$, $m$, $\tau$, $\left\{\sigma_0,\ldots,\sigma_{l-1}\right\}$}\label{algo:PIPELCG}
\begin{algorithmic}[1]
\State $r_{0}:=b-Ax_{0};$ $v_{0}:= r_{0}/{\|r_{0}\|}_2;$ $z_{0}:=v_{0}; ~ g_{0,0}:=1;$
\For {$i=0,\ldots, m+l$}
\State $ z_{i+1}:=\left\{ \begin{matrix}(A-\sigma_{i}I)z_{i}, & i<l \\ Az_{i}, & i \geq l \end{matrix}\right.$ \hfill \# Compute matrix-vector product 
\If {$i\geq l$}  \hfill \# Finalize dot-products ($g_{j,i-l+1}$)
\State $g_{j,i-l+1} := (g_{j,i-l+1}-\sum_{k=i-3l+1}^{j-1}g_{k,j}g_{k,i-l+1})/g_{j,j};$ \qquad $j=i-2l+2,\ldots,i-l$ \hfill \# Update transformation matrix
\State $g_{i-l+1,i-l+1}:= \sqrt{g_{i-l+1,i-l+1}-\sum_{k=i-3l+1}^{i-l}g_{k,i-l+1}^2};$ 
\State \# Check for breakdown and restart if required \hfill \# Square root breakdown check
\If {$i<2l$}
\State $\gamma_{i-l}:=(g_{i-l,i-l+1}+\sigma_{i-l}g_{i-l,i-l} -g_{i-l-1,i-l}\delta_{i-l-1})/g_{i-l,i-l};$ \hfill \# Add column to tridiagonal matrix
\State $\delta_{i-l}:=g_{i-l+1,i-l+1}/g_{i-l,i-l};$
\Else
\State $\gamma_{i-l}:=(g_{i-l,i-l}\gamma_{i-2l}+g_{i-l,i-l+1}\delta_{i-2l} -g_{i-l-1,i-l}\delta_{i-l-1})/g_{i-l,i-l};$
\State $\delta_{i-l}:=(g_{i-l+1,i-l+1}\delta_{i-2l})/g_{i-l,i-l};$
\EndIf
\State \textbf{end if}
\State $v_{i-l+1} := (z_{i-l+1} - \sum_{j=i-3l+1}^{i-l} g_{j,i-l+1} v_{j})/g_{i-l+1,i-l+1};$ \hfill \# Add Krylov subspace basis vector 
\State $z_{i+1} := (z_{i+1} - \gamma_{i-l} z_{i}- \delta_{i-l-1} z_{i-1})/\delta_{i-l};$ \hfill \# Add auxiliary basis vector
\EndIf
\State \textbf{end if}
\State $g_{j,i+1}:=\left\{ \begin{matrix}
	(z_{i+1},v_{j}); & j=\max(0,i-2l+1),\ldots,i-l+1 \\ 
	(z_{i+1},z_{j}); & j=i-l+2,\ldots,i+1 \end{matrix} 
	\right.$ \hfill \text{\# Initiate dot-products ($g_{j,i+1}$)} 
\If {$i=l$} 
\State $\eta_{0}:=\gamma_{0};$ \quad $\zeta_{0}:={\|r_{0}\|}_2;$ \quad $p_{0}:=v_0/\eta_0;$
\Else \, \textbf{if} {$i\geq l+1$} \textbf{then}
\State $\lambda_{i-l}:=\delta_{i-l-1}/\eta_{i-l-1};$  \hfill \text{\# Factorize tridiagonal matrix} 
\State $\eta_{i-l}:=\gamma_{i-l}-\lambda_{i-l}\delta_{i-l-1};$ 
\State $\zeta_{i-l}=-\lambda_{i-l}\zeta_{i-l-1};$  \hfill \text{\# Compute recursive residual norm} 
\State $p_{i-l}=(v_{i-l}-\delta_{i-l-1}p_{i-l-1})/\eta_{i-l};$ \hfill \text{\# Update search direction} 
\State $x_{i-l}=x_{i-l-1}+\zeta_{i-l-1}p_{i-l-1};$ \hfill \text{\# Update approximate solution} 
\If {$|\zeta_{i-l}|/\|r_0\| < \tau$} RETURN; \textbf{end if}  \hfill \# Check convergence criterion
\EndIf 
\EndIf
\State \textbf{end if}
\EndFor 
\end{algorithmic}
}
\end{algorithm*}
Once the basis $V_{i-l+1}$ has been constructed, the solution $x_{i-l}$ can be updated based on a search direction $p_{i-l}$, following the classic derivation of D-Lanczos in \cite{saad2003iterative}, Sec.~6.7.1. The Ritz-Galerkin condition
\begin{align} \label{eq:lanczos}
  0 &= V_m^T r_m  = V_m^T (r_0 - A V_m y_m) \notag \\
		&= V^T_m \left(V_m {\|r_0\|}_2 e_1\right) - T_m y_m = {\|r_0\|}_2 e_1 - T_m y_m,
\end{align}
implies $y_m = T_{m}^{-1} {\|r_0\|}_2 e_1$.
The LU-factorization of the tridiagonal matrix $T_{i-l+1} = L_{i-l+1} U_{i-l+1}$ is given by
\begin{equation} \label{eq:LU}  
\begin{pmatrix} 
1 & & & \\ 
\lambda_{1} & 1 & &  \\  
& \ddots & \ddots &  \\ 
& & \lambda_{i-l} & 1  
\end{pmatrix}
\begin{pmatrix} 
\eta_{0} & \delta_{0} & & \\ 
& \eta_{1} & \ddots & \\  
& & \ddots & \delta_{i-l-1} \\ 
& & & \eta_{i-l} 
\end{pmatrix}.
\end{equation}
Note that $\gamma_{0}=\eta_{0}$. It follows from \eqref{eq:LU} that the elements of the lower/upper triangular matrices $L_{i-l+1}$ and $U_{i-l+1}$ are given by (with $1 \leq j \leq i-l$)
\begin{equation}
\lambda_{j}=\delta_{j-1}/\eta_{j-1} \quad \text{and} \quad \eta_{j}= \gamma_{j}-\lambda_{j}\delta_{j-1}. 
\end{equation}
Expression \eqref{eq:lanczos} indicates that the approximate solution $x_{i-l+1}$ equals 
\begin{align} \label{eq:x_a}
x_{i-l+1} &= x_0 + V_{i-l+1} y_{i-l+1} \notag \\
					&=  x_{0} + V_{i-l+1} U_{i-l+1}^{-1}L_{i-l+1}^{-1}\|r_{0}\|_{2}e_{1} \notag \\
					&= x_{0}+P_{i-l+1} q_{i-l+1},
\end{align}
where $P_{i-l+1} = V_{i-l+1} U_{i-l+1}^{-1}$ and $q_{i-l+1} = L_{i-l+1}^{-1}\|r_{0}\|_{2}e_{1}$. Note that $p_0 = v_0/\eta_0$. The columns $p_j$ (for $1 \leq j \leq i-l$) of $P_{i-l+1}$ can be computed recursively. From $P_{i-l+1} U_{i-l+1} = V_{i-l+1}$ it follows 
\begin{equation}\label{eq:pdir}
p_j = \eta_j^{-1}(v_j-\delta_{j-1} p_{j-1}), \qquad 1 \leq j \leq i-l.
\end{equation}
Denoting the 
vector $q_{i-l+1}$ by $\left[\zeta_0,\ldots,\zeta_{i-l}\right]^T$, it follows from $L_{i-l+1} q_{i-l+1}=\|r_{0}\|_{2}e_{1}$ that $\zeta_0 = {\|r_0\|}_2$ and $\zeta_j = -\lambda_{j}\zeta_{j-1}$ for $1\leq j \leq i-l$. 
Using the search direction $p_{i-l-1}$ and the scalar $\zeta_{i-l-1}$, 
the approximate solution $x_{i-l}$ is updated using the recurrence relation: 
\begin{equation} \label{eq:x_a2}
x_{i-l} =x_{i-l-1} + \zeta_{i-l-1}p_{i-l-1}.
\end{equation}
The above expressions are combined in Alg.\,\ref{algo:PIPELCG}.
Once the initial pipeline for $z_0,\ldots,z_l$ has been filled, the relations \eqref{eq:z_rec}-\eqref{eq:v_rec} are used 
to recursively compute the basis vectors $v_{i-l+1}$ and $z_{i+1}$ in iterations $i \geq l$ (see lines 15-16).
The scalar results of the global reduction phase 
(line 18) are required $l$ iterations later (line 5-6). In every iteration global communication is thus overlapped with the computational work of $l$ subsequent iterations, forming the heart of the communication hiding p($l$)-CG algorithm.

\begin{remark} \label{remark:residual}
\textbf{Residual norm in p($l$)-CG.}
Note that the residual $r_j = b-Ax_j$ is not computed in Alg.\,\ref{algo:PIPELCG}, but its norm is characterized by the quantity $|\zeta_j| = \|r_j\|$ for $0 \leq j \leq i-l$. 
This quantity can be used to formulate a stopping criterion for the p($l$)-CG iteration, see Alg.\,\ref{algo:PIPELCG} line 27.
\end{remark}

\begin{remark} \label{remark:dotpr}
\textbf{Dot products in p($l$)-CG.}
Although Alg.\,\ref{algo:PIPELCG}, line 18 indicates that in each iteration $i \geq (2l+1)$ a total of $(2l+1)$ dot products need to be computed, the number of dot product computations can be limited to $(l+1)$ by exploiting the symmetry of the matrix $G_{i+1}$, see expression \eqref{eq:symmetry_G}. Since $g_{j,i+1} = g_{i-l+1,j+l}$ for $j \leq i+1$, only the dot products $(z_{j},z_{i+1})$ for $j=i-l+2,\ldots,i+1$
and the $l$-th upper diagonal element $(v_{i-l+1},z_{i+1})$ need to be computed in iteration $i$.
\end{remark}

\subsection{On the conditioning of the auxiliary basis 
} \label{sec:conditioning}

As $V_i$ is an orthonormal basis, the transformation $Z_i = V_i G_i$ can be interpreted as a QR factorization of the auxiliary basis $Z_i$. Moreover, it holds that $G_i^T G_i$ is the Cholesky factorization of $Z_i^T Z_i$, since
\begin{equation} \label{eq:ZTZ}
	Z_i^T Z_i = G_i^T V_i^T V_i G_i = G_i^T G_i.
\end{equation}
The elements of the transformation matrix $G_i$ are computed on lines 5-6 of Alg.\,\ref{algo:PIPELCG} precisely by means of this Cholesky factorization. This observation leads to the following two essential insights related to the conditioning of the basis $Z_i$.

\begin{remark} \label{remark1}
\textbf{Square root breakdowns in p($l$)-CG.} The auxiliary basis vectors $z_j$ are defined as $P_l(A) v_{j-l}$, but the basis $Z_{i}$ is in general not orthogonal. Hence, vectors $z_j \in Z_{i}$ are not necessarily linearly independent. In particular for longer $l$, different $z_j$ vectors are expected to become more and more aligned. This leads to $Z_i^T Z_i$ being ill-conditioned, approaching singularity as $i$ increases. 
The effect is the most pronounced when $\sigma_0 = \ldots = \sigma_{l-1} = 0$, in which case $P_l(A) = A^l$. Shifts $\sigma_j$ can be set to improve the conditioning of $Z_i^T Z_i$, see also Remark \ref{remark2}.

When for certain $i$ the matrix $Z_i^T Z_i$ becomes (numerically) singular, the Cholesky factorization procedure in p($l$)-CG will fail. This may manifest in the form of a square root breakdown on line 7 in Alg.\,\ref{algo:PIPELCG} when the root argument $g_{i-l+1,i-l+1}-\sum_{k=i-3l+1}^{i-l}g_{k,i-l+1}^2$ becomes negative. Numerical round-off errors may increase the occurrence of these breakdowns in practice. When a breakdown occurs in p($l$)-CG the iteration is restarted, in analogy to the GMRES algorithm, although it should be noted that the nature of the breakdown in both algorithms is quite different. Evidently, the restarting strategy may delay convergence compared to standard CG, where no square root breakdowns occur.
\end{remark}

\begin{remark} \label{remark2}
\noindent \textbf{Choice of the auxiliary basis 
and relation to the shifts. 
} It follows from the Cholesky factorization \eqref{eq:ZTZ} that the inverse of the transformation matrix $G_i$ is
\begin{equation}
	G_i^{-1} = (Z_i^T Z_i)^{-1} G_i^T.
\end{equation}
The conditioning of the matrix $G_i^{-1}$ is thus determined by the conditioning of $Z_i^T Z_i$. Furthermore, it holds that
\begin{equation}
	Z_{l+1:i}^T Z_{l+1:i} = (P_l(A) V_{i-l})^T P_l(A) V_{i-l} = V_{i-l}^T \, P_l(A)^2 \, V_{i-l},
\end{equation}
where $Z_{l+1:i}$ is a part of the basis $Z_i$ obtained by dropping the first $l$ columns.
Hence, the polynomial $P_l(A)^2$ has a major impact on the conditioning of the matrix $G_i^{-1}$, which in turn plays a crucial role in the propagation of local rounding errors in the p($l$)-CG algorithm \cite{cornelis2017communication}, see Section \ref{sec:local} of the current work. This observation indicates intuitively why $\|P_l(A)\|_2$ should preferably be as small as possible, which can be achieved by choosing appropriate values for the shifts $\sigma_j$. Optimal shift values in the sense of minimizing the Euclidean $2$-norm of $P_l(A)$ are the Chebyshev shifts \cite{hoemmen2010communication,ghysels2013hiding,cornelis2017communication} (for $i=0,\ldots,l-1$):
\begin{equation} \label{eq:chebyshev}
  \sigma_{i} = \frac{\lambda_{\max}+\lambda_{\min}}{2}+\frac{\lambda_{\max}-\lambda_{\min}}{2} \cos\left(\frac{(2i+1)\pi}{2l}\right),
\end{equation}
which are used throughout this work. This requires a notion of the largest (
smallest) eigenvalue $\lambda_{\max}$ (resp.~$\lambda_{\min}$), which can be estimated a priori, e.g.~by 
a few 
power method iterations.
\end{remark}

\subsection{
Loss of orthogonality and attainable accuracy
}

This work is motivated by the observation that two main issues affect the convergence of pipelined CG methods: loss of basis vector orthogonality and inexact Lanczos relations. We comment briefly on both issues from a high-level point of view and clearly mark the scope of this work. Important insights about the similarities and differences between classic CG and p($l$)-CG are highlighted before going into more details on the numerics in Sections \ref{sec:analyzing} and \ref{sec:stabilizing}.

\subsubsection{Loss of orthogonality} It is well-known that in finite precision arithmetic the orthogonality of the Krylov subspace basis $V_i$, i.e.~$V^T_i V_i = I_i$ (identity matrix) may not hold exactly. Inexact orthogonality may appear in every variant of the CG algorithm, see \cite{liesen2012krylov}, in particular in the D-Lanczos\footnote{\textbf{Note:} The D-Lanczos (short for ``\emph{direct Lanczos}'') algorithm is a variant of the CG method that is equivalent to the latter in exact arithmetic, save for the solution of the system $T_i y_i = \|r_0\| e_1$ which is computed by using Gaussian elimination in D-Lanczos. The D-Lanczos method is the basic Krylov subspace method from which the p($l$)-CG method was derived, see 
\cite{cornelis2017communication}, Section 2, for details.} algorithm \cite{saad2003iterative}, where a new basis vector is constructed by orthogonalizing with respect to the previous two basis vectors, as well as in the related p($l$)-CG method, Alg.\,\ref{algo:PIPELCG}. Loss of orthogonality typically leads to \emph{delay of convergence}, 
meaning the residual deviates from the one in the scenario in which orthogonality would not be lost.

We use a \emph{notation with bars} to designate variables that are actually computed in a finite precision implementation of the algorithm.
 The key relation for the Conjugate Gradient method is the Ritz-Galerkin condition:
\begin{align}
	\bar{V}_i^T (b-A\bar{x}_i) &= \bar{V}_i^T (\bar{r}_0 - A \bar{V}_i \bar{T}_i^{-1} \|\bar{r}_0\| e_1) \notag \\
														 &= \bar{V}_i^T \bar{v}_0 \|\bar{r}_0\| - \bar{V}_i^T A \bar{V}_i \bar{T}_i^{-1} \|\bar{r}_0\| e_1 = 0.
\end{align}
This equality only holds under the assumption that $\bar{V}_{i}^T A \bar{V}_{i}  = \bar{T}_{i}$ which requires $\bar{V}_i^T \bar{V}_{i+1} = I_{i,i+1}$. Note that in finite precision arithmetic the convergence delay can be observed in both the actual residual norm $\|b-A\bar{x}_i\|$ as well as the recursively computed residual norm $\|\bar{r}_i\|$, since both quantities are based on the (possibly non-orthogonal) basis $\bar{V}_{i+1}$, see Fig.\,\ref{fig:figure1a}-\ref{fig:figure2a} (discussed further in Section \ref{sec:scope}).

\subsubsection{The inexact Lanczos relation} The basis vectors in the pipelined CG algorithm are not computed explicitly using the Lanczos relation \eqref{eq:Arnoldi_v}. Rather, they are computed recursively, see \eqref{eq:v_rec}, to avoid the computation of additional \textsc{spmv}s. In finite precision, local rounding errors in the recurrence relation may contaminate the basis $\bar{V}_i$, such that the Lanczos relation $A\bar{V}_i - \bar{V}_{i+1}\bar{T}_{i+1,i} = 0$ is no longer valid. Moreover, due to propagation of local rounding errors, $\|A\bar{V}_i - \bar{V}_{i+1}\bar{T}_{i+1,i}\|$ may grow dramatically as the iteration proceeds. 
Using the classic model for floating point arithmetic with machine precision $\epsilon$ \cite{paige1976error,greenbaum1992predicting,greenbaum1997estimating,wilkinson1994rounding}, the round-off error on basic computations on the matrix $A \in \mathbb{R}^{n \times n}$, vectors $v$, $w$ and a scalar $\alpha$ are bounded by
\begin{align*}
	\| \alpha v - \text{fl}(\alpha v) \| &\leq \| \alpha v \| \, \epsilon =  |\alpha| \, \|v\| \, \epsilon, \\
	\| v + w - \text{fl}(v + w) \| &\leq (\|v\| + \|w\|) \, \epsilon,\\
	\| Av - \text{fl}(Av) \| &\leq \mu\sqrt{n} \, \|A\| \, \|v\| \, \epsilon, 
\end{align*}
Here $\text{fl}(\cdot)$ indicates the finite precision floating point representation, $\mu$ is the maximum number of nonzeros in any row of $A$, and the norm $\|\cdot\|$ represents the Euclidean 2-norm. Under this model the recurrence relations for $\bar{x}_i$ and $\bar{p}_i$ in a finite precision implementation of p($l$)-CG are
\begin{align} 
	\bar{x}_i &= \bar{x}_{i-1} + \bar{\zeta}_{i-1} \bar{p}_{i-1} + \xi^{\bar{x}}_i = \bar{x}_0 + \bar{P}_{i} \bar{q}_i + \Xi_i^{\bar{x}} \, \bold{1}, \label{eq:exp_x}\\
	\bar{p}_i &= (\bar{v}_{i} - \bar{\delta}_{i-1} \bar{p}_{i-1})/\bar{\eta}_i + \xi^{\bar{p}}_i, \label{eq:exp_p}
\end{align}
with expression \eqref{eq:exp_p} translating in matrix notation to
\begin{equation}
	\bar{V}_{i} = \bar{P}_{i} \bar{U}_{i} + \Xi^{\bar{p}}_{i}. \label{eq:exp_p2}
\end{equation}
Recall that in exact arithmetic $P_i = [p_0, \ldots, p_{i-1}] = V_i U^{-1}_i$.
In these expressions $\Xi_i^{\bar{x}} = [\xi^{\bar{x}}_1,
\ldots,\xi^{\bar{x}}_{i}]$ and $ \Xi^{\bar{p}}_{i} = -[\bar{\eta}_0\xi^{\bar{p}}_0,
,\ldots,\bar{\eta}_{i-1}\xi^{\bar{p}}_{i-1}]$ are local rounding errors which are bounded by $\|\xi^{\bar{x}}_i\| \leq (\|\bar{x}_{i-1}\| + 2 \, |\bar{\zeta}_{i-1}| \, \|\bar{p}_{i-1}\|) \, \epsilon$ and $\|\xi^{\bar{p}}_i\| \leq (2 / \bar{\eta}_i \, \|\bar{v}_{i-1}\|  + 3 \, |\bar{\delta}_{i-1}|/ \bar{\eta}_i  \, \|\bar{p}_{i-1}\|) \, \epsilon$, and $\bold{1} = [1,1,\ldots,1]^T$.
The actual residual satisfies the following relations in a finite precision setting: 
\begin{align} \label{eq:resgap}
	b-A\bar{x}_i 
			&= \bar{r}_0 - A \bar{V}_i \bar{U}_i^{-1} \bar{q}_i + \overbrace{A \Xi_i^{\bar{p}} \bar{U}^{-1}_{i} \bar{q}_i - A \Xi_i^{\bar{x}} \, \bold{1} + \xi^{\bar{r}}_0}^{\text{Local Rounding Errors (LRE)}} \notag \\
			&= \bar{r}_0 - \bar{V}_{i+1} \bar{T}_{i+1,i} \bar{U}_i^{-1} \bar{q}_i \notag \\
					& \qquad \qquad - (A\bar{V}_i - \bar{V}_{i+1}\bar{T}_{i+1,i}) \bar{U}_i^{-1} \bar{q}_i + \text{LRE} \notag \\
			&= \bar{r}_i - (A\bar{V}_i - \bar{V}_{i+1}\bar{T}_{i+1,i}) \bar{U}_i^{-1} \bar{q}_i + \text{LRE}
\end{align}
where $\xi^{\bar{r}}_0 = (b-A\bar{x}_0) - \bar{r}_0$. The recursively computed residual $\bar{r}_i = -\bar{\delta}_{i-1} (e_i^T \bar{U}_i^{-1} \bar{q}_i) \bar{v}_i = \bar{\xi}_i \bar{v}_i $ that appears in expression \eqref{eq:resgap} tends to zero.
The actual residual norm $\|b-A\bar{x}_i\|$, on the other hand, stagnates around $\|(A\bar{V}_i - \bar{V}_{i+1}\bar{T}_{i+1,i}) \bar{U}_i^{-1} \bar{q}_i\|$, a quantity referred to as the \emph{maximal attainable accuracy} of the method. 
The difference between the norm of the actual residual and the recursively computed residual is illustrated in Fig.\,\ref{fig:figure1a}.
A detailed analysis of the deviation from the Lanczos relation in finite precision (``\emph{inexact Lanczos}'') can be found in Section \ref{sec:analyzing}.

\subsubsection{Scope and limitations of this manuscript}\label{sec:scope}   
\begin{figure*}
    \centering
    \begin{subfigure}[b]{0.49\textwidth}
				\includegraphics[width=\textwidth]{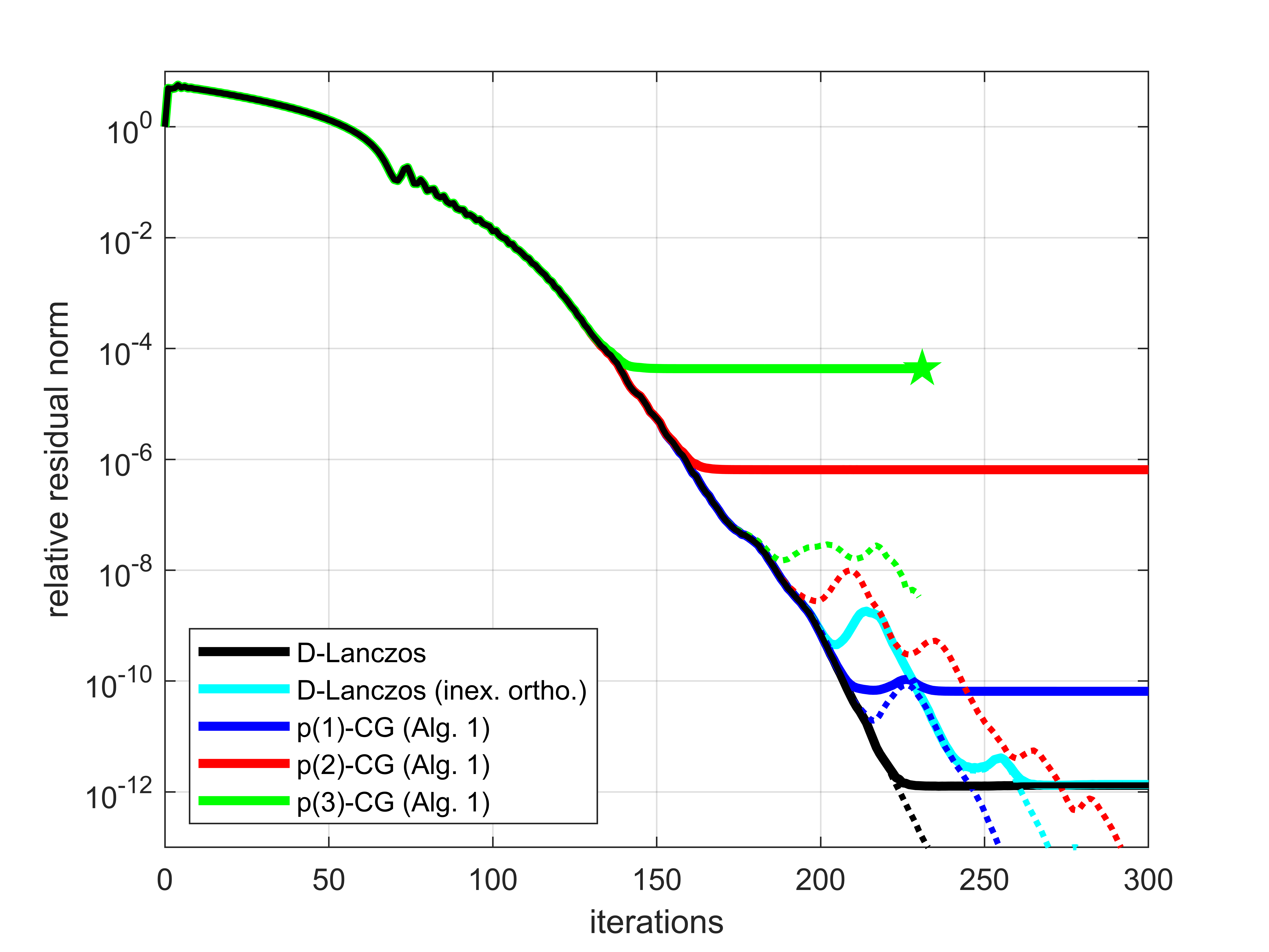} 
        \caption{}
        \label{fig:figure1a}
    \end{subfigure}
    \begin{subfigure}[b]{0.49\textwidth}
				\includegraphics[width=\textwidth]{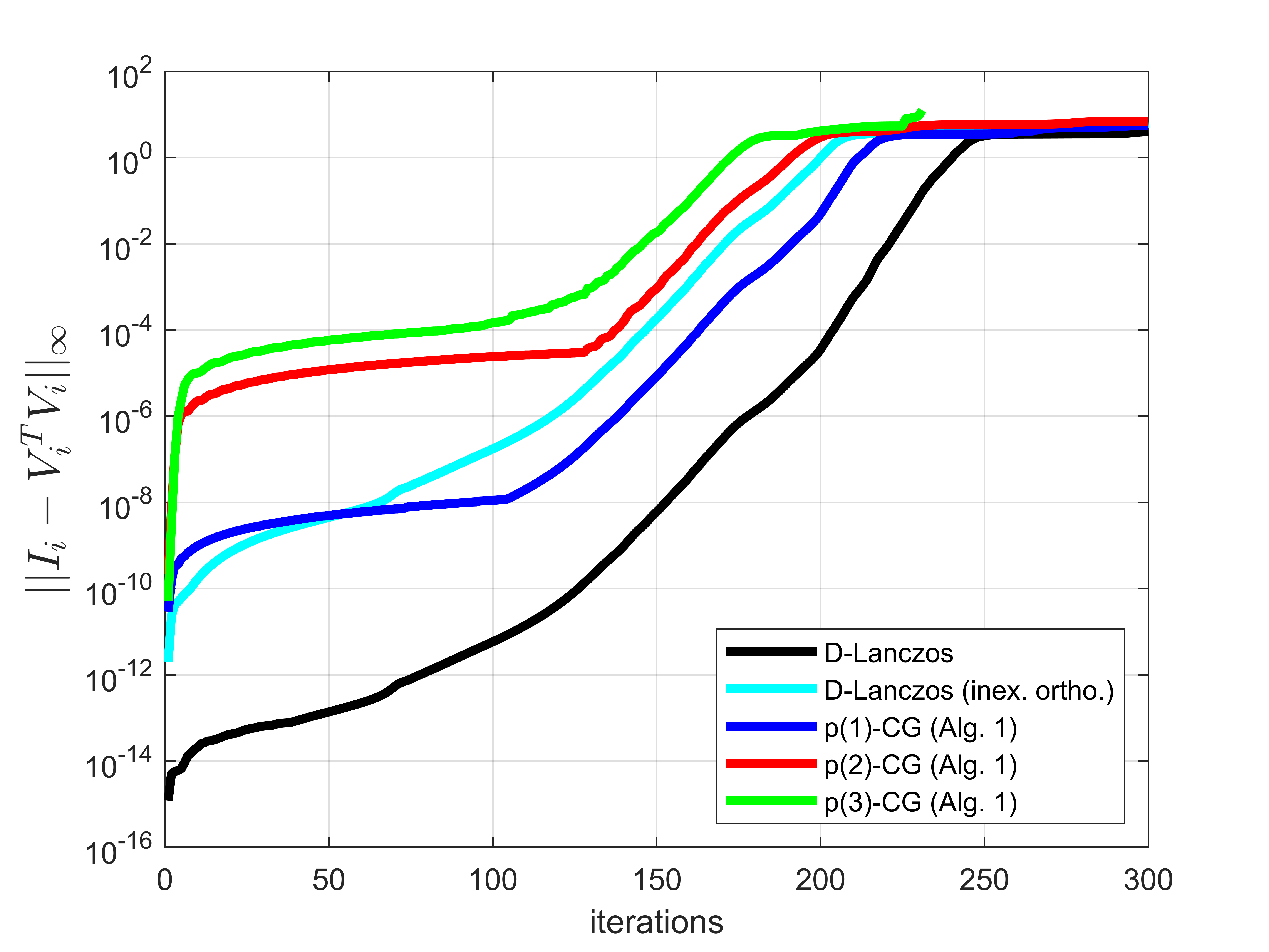}
        \caption{}
        \label{fig:figure1b}
    \end{subfigure}
		\vspace{-0.2cm}
\caption{
\small
Reduced attainable accuracy and loss of orthogonality in CG variants.  
(a) 
Actual residual norm $\|b-A\bar{x}_i\|/\|b\|$ (full line) and recursive residual norm $\|\bar{r}_i\|/\|b\|$ (dotted line) 
for a 2D $100 \times 100$ unknowns Laplace problem. 
Square root breakdown 
indicated by $\bigstar$, 
see Remark \ref{remark1}. 
(b) Basis orthogonality characterized by $\|I - \bar{V}_i^T \bar{V}_i\|_{\infty}$.
For D-Lanczos (cyan) 
orthogonality was modified by setting $\bar{\gamma}_i = (A\bar{v}_i,\bar{v}_i) (1 + \varepsilon)$ with $\varepsilon = 1.0$e-$11$ in each iteration.
}
\label{fig:figure1}
\end{figure*}

The issue of inexact Lanczos relations in p($l$)-CG  
is timely and deserving of attention. Loss of accuracy resulting from the inexact Lanczos relation 
has long been a limiting factor in applying p($l$)-CG 
and related algorithms in practice. 
Fig.\,\ref{fig:figure1a} illustrates how the norms of the actual residuals $\|b-A\bar{x}_i\|$ stagnate while the recursively computed residual norms $\|\bar{r}_i\|$ continue to decrease. For p($l$)-CG local rounding errors in the recurrence relations are propagated leading to reduced attainable accuracy compared to D-Lanczos.
Moreover, while 
loss of orthogonality also warrants further investigation, 
this issue is not exclusive to pipelined methods. Delayed convergence is observed in classic Krylov subspace methods also, see Fig.\,\ref{fig:figure1}, whereas loss of attainable accuracy 
is not. The issue 
could be addressed by e.g.~re-orthogonalizing the basis, see \cite{swirydowicz2018low}. However, communication-reducing methods are not particularly suitable to include re-orthogonalization, since this introduces additional global reduction phases. 

Although loss of orthogonality does not originate 
from applying the pipelining technique, it may be more pronounced for pipelined 
methods compared to their classic counterparts, see Fig.\,\ref{fig:figure1b}. However, Fig.\,\ref{fig:figure1a} indicates that the effect of the inexact Lanczos relation on convergence 
is much more dramatic than the effect of inexact orthogonality for all pipeline lengths $l$. 
 This manuscript thus focuses on improving the numerical stability of the p($l$)-CG method by neutralizing the propagation of local rounding errors in the recursively computed basis vector updates. As such, this work 
proposes a key step towards 
a numerically stable communication hiding variant of the 
CG method. 

\section{Analyzing rounding error propagation 
} \label{sec:analyzing} 

This section recaps the analysis of local rounding errors that stem from the recurrence relations in the pipelined p($l$)-CG method, Alg.\,\ref{algo:PIPELCG}. It aims to precisely explain the source of the loss of accuracy observed for the p($l$)-CG method. 
The methodology for the analysis is similar to the one used in classic works by Paige \cite{
paige1976error,paige1980accuracy}, Greenbaum \cite{greenbaum1989behavior,greenbaum1992predicting,greenbaum1997estimating}, Gutknecht \cite{gutknecht2000accuracy}, Strakos \cite{
strakovs2005error}, Meurant \cite{meurant2006lanczos},  Sleijpen \cite{sleijpen1996reliable,sleijpen2001differences}, Van der Vorst \cite{van2000residual}, Higham \cite{higham2002accuracy}, and others.

Finite precision variants of the exact scalar and vector variables introduced in Section \ref{sec:pipelcg} are denoted by a bar symbol in this section.
We differentiate between ``\emph{actual}'' vector variables, which satisfy the Lanczos relations exactly but are not computed in the algorithm, and ``\emph{recursively computed}'' variables, which contain machine-precision sized round-off errors related to finite precision computations.

\subsection{Local rounding error behavior in finite precision 
} \label{sec:local}

For any $j \geq 0$ the \emph{true basis vector}, denoted by $\bar{\bold{v}}_{j+1}$, satisfies the Lanczos relation \eqref{eq:v_ex} exactly, that is, for $0 \leq j < i-l$:
\begin{equation} \label{eq:v_arnoldi}
  \bar{\bold{v}}_{j+1}= ( A\bar{v}_{j} - \bar{\gamma}_j \bar{v}_j - \bar{\delta}_{j-1} \bar{v}_{j-1} ) / \bar{\delta}_j,
\end{equation}
without the addition of a local rounding error. This vector is not actually computed in the p($l$)-CG algorithm.
Instead, the \emph{computed basis vector} $\bar{v}_{j+1}$ is calculated from the finite precision variant of relation 
\eqref{eq:v_rec} for $0 \leq j < i-l$, i.e.:
\begin{equation} \label{eq:vbar_rec}
  \bar{v}_{j+1} = \left( \bar{z}_{j+1} - \sum_{k=j-2l+1}^{j} \bar{g}_{k,j+1} \bar{v}_{k} \right) / \bar{g}_{j+1,j+1} + \xi^{\bar{v}}_{j+1},
\end{equation}
where the size of the local rounding errors 
is bounded in terms of machine precision: $\| \xi^{\bar{v}}_{j+1} \| \leq ( 2 \, \|\bar{z}_{j+1}\|/|\bar{g}_{j+1,j+1}| + 3 \sum_{k=j-2l+1}^{j} |\bar{g}_{k,j+1}|/|\bar{g}_{j+1,j+1}| \, \|\bar{v}_k\| ) \epsilon$. 
Let $\bar{V}_{j+1} = [\bar{v}_0, \bar{v}_1, \ldots, \bar{v}_j]$ and $\bar{\bold{V}}_{j+1} = [\bar{\bold{v}}_0, \bar{\bold{v}}_1, \ldots, \bar{\bold{v}}_j]$.
Relation \eqref{eq:v_arnoldi} alternatively translates to the following formulation in matrix notation (with $1 \leq j \leq i-l$):
\begin{equation} \label{eq:AV_BAR}
  A \bar{V}_j = \bar{V}_{j+1} \bar{T}_{j+1,j} + (\bar{\bold{V}}_{j+1}-\bar{V}_{j+1}) \bar{\Delta}_{j+1,j}, 
\end{equation}
where $\bar{\Delta}_{j+1,j}$ is a $(j+1)$-by-$j$ rectangular matrix holding the entries $\{\bar{\delta}_0, \ldots, \bar{\delta}_{j-1}\}$ directly below the main diagonal.
The matrix $(\bar{\bold{V}}_{j+1}-\bar{V}_{j+1}) = [\bar{\bold{v}}_0 - \bar{v}_0, \bar{\bold{v}}_1 - \bar{v}_1, \ldots, \bar{\bold{v}}_j - \bar{v}_j]$ collects the \emph{gaps} between the actual and recursively computed basis vectors, which quantify the deviation from the Lanczos relation in the finite precision setting. These gaps are crucial in describing the propagation of local rounding errors throughout the p($l$)-CG algorithm and are directly linked to the gap between the actual and recursively computed residuals, see expression \eqref{eq:resgap}. 
From \eqref{eq:vbar_rec} one obtains 
\begin{equation} \label{eq:Z_BAR}
  \bar{Z}_j = \bar{V}_j \bar{G}_j + \Xi^{\bar{v}}_j, \qquad 1 \leq j \leq i-l,
\end{equation}
where $\Xi^{\bar{v}}_j = [0, \, \bar{g}_{1,1}\xi^{\bar{v}}_1, \, \ldots, \, \bar{g}_{j-1,j-1}\xi^{\bar{v}}_{j-1}]$ collects the local rounding errors. The computed auxiliary vector $\bar{z}_{j+1}$ satisfies a finite precision version of the recurrence relation \eqref{eq:z_rec}, which summarizes to
\begin{equation} \label{eq:zbar_rec}
  \bar{z}_{j+1} = \left\{
	\begin{matrix}
			(A-\sigma_jI) \, \bar{z}_j + \xi^{\bar{z}}_{j+1}, \hfill  0 \leq j < l, \vspace{0.2cm}\\
			\displaystyle{\frac{ A\bar{z}_{j} - \bar{\gamma}_{j-l} \bar{z}_j - \bar{\delta}_{j-l-1} \bar{z}_{j-1} }{ \bar{\delta}_{j-l} }} + \xi^{\bar{z}}_{j+1}, \quad l \leq j < i, 
	\end{matrix} \right.
\end{equation}
where $\xi^{\bar{z}}_{j+1}$ is the local rounding error which can again be bounded in terms of machine precision $\epsilon$. Expression \eqref{eq:zbar_rec} can be formulated in matrix notation as:
\begin{equation} \label{eq:AZ_BAR}
  A \bar{Z}_j = \bar{Z}_{j+1} \bar{B}_{j+1,j} + \Xi^{\bar{z}}_j, \qquad 1 \leq j \leq i.
\end{equation}
Furthermore, the following matrix relations hold between the scalar coefficients $\bar{\gamma}_j$ and $\bar{\delta}_j$ in Alg.\,\ref{algo:PIPELCG}:
\begin{equation} \label{eq:rec_coeff}
  \bar{G}_{j+1}  \bar{B}_{j+1,j} = \bar{T}_{j+1,j} \bar{G}_j, \qquad 1 \leq j \leq i-l.
\end{equation}
Subsequently, using expressions \eqref{eq:AV_BAR}, \eqref{eq:Z_BAR}, \eqref{eq:AZ_BAR} and \eqref{eq:rec_coeff} and it is derived that the gaps on the basis vectors are given by
\begin{equation} \label{eq:gap_plcg2}
\bar{\bold{V}}_{j+1} - \bar{V}_{j+1} 
				= (\Xi^{\bar{z}}_j - A \Xi^{\bar{v}}_j + \Xi^{\bar{v}}_{j+1} \bar{B}_{j+1,j} ) \, \bar{G}^{-1}_j  \bar{\Delta}^{+}_{j+1,j},
\end{equation}
where $\bar{\Delta}_{j+1,j}^{+} = (\bar{\Delta}_{j+1,j}^* \bar{\Delta}_{j+1,j})^{-1} \bar{\Delta}_{j+1,j}^*$ should be interpreted as a Moore-Penrose (left) pseudo-inverse. Hence, the local rounding errors in this expression are possibly amplified by the entries of the matrix $\bar{G}^{-1}_j\bar{\Delta}^{+}_{j+1,j}$, which may lead to loss of attainable accuracy for the p($l$)-CG method. The inexact Lanczos relation may in turn give rise to a growing gap between the computed and actual residual on the solution, see \eqref{eq:resgap}. It is clear from expression \eqref{eq:gap_plcg2} that the conditioning of the matrix $\bar{G}^{-1}_j$ plays a crucial role in the rounding error propagation in the p($l$)-CG algorithm as indicated in Section \ref{sec:conditioning}, see Remark \ref{remark2}. 



\subsection{Toward stability by using the Lanczos relation 
} \label{sec:countermeasures}

Section \ref{sec:local} shows that the recurrence relation \eqref{eq:vbar_rec} is the main cause for the amplification of local rounding errors throughout the p($l$)-CG algorithm. Moreover, the possibly ill-conditioned matrix $\bar{G}_j$ that is used construct the basis $\bar{V}_j$, see expression \eqref{eq:Z_BAR}, may be detrimental for convergence. A straightforward countermeasure would be to eliminate $\bar{G}_j$ in the construction of the basis. This can be achieved by simply replacing the recurrence relation \eqref{eq:vbar_rec} by the original Lanczos relation, i.e., for $0 \leq j < i - l$:
\begin{equation} \label{eq:vbar_stab}
	\bar{v}_{j+1} = ( A\bar{v}_{j} - \bar{\gamma}_j \bar{v}_j - \bar{\delta}_{j-1} \bar{v}_{j-1} ) / \bar{\delta}_j + \psi^{\bar{v}}_{j+1}.
\end{equation} 
Here $\psi^{\bar{v}}_{j+1}$ represents a local rounding error which is generally different from the error $\xi^{\bar{v}}_{j+1}$ occurring in expression \eqref{eq:vbar_rec}.
Recurrence relation \eqref{eq:vbar_stab} 
can alternatively be written as 
\begin{equation}
	A\bar{V}_j = \bar{V}_{j+1} \bar{T}_{j+1,j} - \Psi^{\bar{v}}_{j+1} \bar{\Delta}_{j+1,j},
\end{equation}
with $1 \leq j \leq i-l$.
The gap between the true basis vector $\bar{\bold{v}}_{j+1}$ and the computed basis vector $\bar{v}_{j+1}$ then reduces to
\begin{equation} \label{eq:gap_stab}
	\bar{\bold{V}}_{j+1} - \bar{V}_{j+1} = - \Psi^{\bar{v}}_{j+1}. 
\end{equation}
By using recurrence \eqref{eq:vbar_stab} for $\bar{v}_{j+1}$ instead of relation \eqref{eq:vbar_rec} in Alg.\,\ref{algo:PIPELCG}, no amplification of local rounding errors occurs, see \eqref{eq:gap_stab}, and the influence of rounding errors on attainable accuracy remains limited. However, to use the recurrence relation \eqref{eq:vbar_stab} an additional \textsc{spmv}, i.e.~$A\bar{v}_{j}$, is computed in each iteration of the 
algorithm, leading to an undesirable increase in computational cost. Although the use of expression \eqref{eq:vbar_stab} would not hinder the ability to overlap the global reduction phase with computations (
for pipeline length $l$ the global reduction would simply be overlapped with $2l$ \textsc{spmv}s), we aim to avoid adding \textsc{spmv} computations to the algorithm.

The technique proposed by expression \eqref{eq:vbar_stab} shows similarity to the concept of \emph{residual replacement}, 
which was suggested by several authors as a countermeasure to local rounding error propagation in various multi-term recurrence variants of CG \cite{van2000residual,carson2014residual,cools2018analyzing}. While the idea is valuable, it cannot be implemented in the p($l$)-CG method in its current form, i.e.~using expression \eqref{eq:vbar_stab}, due to the significantly augmented computational cost caused by the additional \textsc{spmv} in each iteration.

\section{Deriving stable recurrence relations 
} \label{sec:stabilizing}

We now present the core technique for obtaining a numerically stable variant of the recurrence relations used in the p($l$)-CG algorithm by introducing additional auxiliary bases and corresponding recurrence relations. Sections \ref{sec:derivation}-\ref{sec:preconditioning} are again written in the exact arithmetic framework in order to derive the algorithm, interluded by a short discussion on computational costs and storage requirements in Section \ref{sec:computational}. We return to the finite precision framework for the analysis of the improved method in Section \ref{sec:rounding}.

\subsection{Derivation of a stable pipelined CG method} \label{sec:derivation}

We introduce a total of $l+1$ bases, denoted by $Z^{(k)}_{i+1}$, where the upper index `$(k)$' (with $0 \leq k \leq l$) labels the different bases and the lower index `$i+1$' indicates the iteration like before.
The basis $Z^{(0)}_{i+1}$ will denote the original Krylov subspace basis, that is: $Z^{(0)}_{i+1} = V_{i+1}$.
 The auxiliary basis vectors $Z^{(l)}_{i+1} =[z^{(l)}_{0},z^{(l)}_{1},\ldots,z^{(l)}_{i}]$ are defined identically to the basis $Z_{i+1}$ in p($l$)-CG, cf.~\eqref{eq:z_ex}, i.e.~$Z^{(l)}_{i+1} = Z_{i+1}$.
In addition, we also define $l-1$ intermediary bases $Z^{(1)}_{i+1}, \ldots, Z^{(l-1)}_{i+1}$ that will enable us to use a variant of the Lanczos relation \eqref{eq:vbar_stab} to recursively update $v_j$, but without the necessity to compute the \textsc{spmv} $Av_{j}$. The auxiliary bases are defined as follows:
\begin{equation} \label{eq:def_zjk}
z^{(k)}_{j} = 
\left\{ 
	\begin{matrix} 
		v_{0}, & j=0, \\ 
		P_{j}(A)v_{0}, & 0<j\leq k, \\ 
		P_{k}(A)v_{j-k}, & k<j\leq i, 
	\end{matrix} 
\right. \quad \text{for}~~ 0 \leq k \leq l,
\end{equation}
where the polynomial is defined by \eqref{eq:poly}.
Note that the first $k+1$ vectors in basis $Z^{(k)}_{j}$ $(j \geq k)$ are identical to the first $k+1$ vectors in all bases $Z^{(k')}_{j}$ with $k' \geq k$.
By definition \eqref{eq:def_zjk} the `zero-th' basis $Z^{(0)}_{j}$ is simply the original basis $V_j$, whereas the $l$-th basis $Z^{(l)}_{j}$ is the auxiliary basis $Z_j$ from the p($l$)-CG method, see Section \ref{sec:pipelcg}.

A crucial relation connects each pair of bases $Z^{(k)}_j$ and $Z^{(k+1)}_j$ (for $j > k$). It holds that
\begin{align} 
  z^{(k+1)}_{j+1} &= P_{k+1}(A) v_{j-k} = (A-\sigma_{k}) P_{k}(A) v_{j-k} \notag \\
									&= (A-\sigma_{k}) z^{(k)}_{j}, \quad j \geq k, \quad 0 \leq k \leq l-1,
\end{align}
which translates into
\begin{equation} \label{eq:all_z_relation}
  Az^{(k)}_j = z^{(k+1)}_{j+1} + \sigma_{k} z^{(k)}_{j}, \quad j \geq k, \quad 0 \leq k \leq l-1.
\end{equation}
By multiplying the original Lanczos relation \eqref{eq:v_ex} for $v_j$ on both sides by the respective polynomial $P_k(A)$ with $1 \leq k \leq l$ and by exploiting the associativity of $A$ and $P_k(A)$, it is straightforward to derive that each auxiliary basis $Z^{(k)}_j$ with $0 \leq k \leq l$ satisfies a Lanczos type recurrence relation:
\begin{equation} \label{eq:all_z_rec}
z^{(k)}_{j+1} = (A z^{(k)}_{j} - \gamma_{j-k} z^{(k)}_{j} - \delta_{j-k-1} z^{(k)}_{j-1}) / \delta_{j-k}, 
\end{equation}
for $j \geq k$ and $0 \leq k \leq l$.
Note that when $k = 0$ expression \eqref{eq:all_z_rec} yields the Lanczos relation \eqref{eq:v_ex} for $v_{j+1}$, whereas setting $k = l$ results in the recurrence relation \eqref{eq:z_rec} for $z_{j+1}$. 

The recursive expressions \eqref{eq:all_z_rec} for the bases $Z^{(0)}_j, \ldots, Z^{(l-1)}_j$ are not particularly useful in practice, since each recurrence relation requires to compute an additional \textsc{spmv} to form the next basis vector. However, using relation \eqref{eq:all_z_relation} the recurrence relations \eqref{eq:all_z_rec} can alternatively be written as:
\begin{equation} \label{eq:all_z_rec2}
z^{(k)}_{j+1} = (z^{(k+1)}_{j+1} + (\sigma_{k} - \gamma_{j-k}) z^{(k)}_{j} - \delta_{j-k-1} z^{(k)}_{j-1}) / \delta_{j-k},
\end{equation}
with $j \geq k$ and $0 \leq k < l$. We stress that only for $k = l$, i.e.~to compute the vectors in the auxiliary basis $Z^{(l)}_j = Z_j$, we use the recursive update 
given by expression \eqref{eq:all_z_rec}:
\begin{equation} \label{eq:specific_z_rec}
z^{(l)}_{j+1} = (A z^{(l)}_{j} - \gamma_{j-l} z^{(l)}_{j} - \delta_{j-l-1} z^{(l)}_{j-1}) / \delta_{j-l}, 
\end{equation}
for $j \geq l$,
which reduces to the recurrence relation \eqref{eq:z_rec}. The recurrence relations \eqref{eq:all_z_rec2} allow us to compute the vector updates for the bases $Z^{(0)}_j, \ldots, Z^{(l-1)}_j$ without the need to compute any additional \textsc{spmv}. Adding the recurrence relations \eqref{eq:all_z_rec2} for the auxiliary bases $Z^{(0)}_j, \ldots, Z^{(l-1)}_j$ to the p($l$)-CG method leads to the stable p($l$)-CG method, Alg.\,\ref{algo:PIPELCGSTAB}.

\begin{algorithm*}[t]
{\small
\caption{Stable $l$-length pipelined p($l$)-CG \hfill 
\textbf{Input:} $A$, $b$, $x_0$, $l$, $m$, $\tau$, $\left\{\sigma_0,\ldots,\sigma_{l-1}\right\}$}\label{algo:PIPELCGSTAB}
\begin{algorithmic}[1]
\State $r_{0}:=b-Ax_{0};$ $z^{(l)}_{0}:= z^{(l-1)}_{0}:= \,\ldots\, := z^{(1)}_{0}:=z^{(0)}_{0} := r_{0}/{\|r_{0}\|}_2; ~ g_{0,0}:=1;$
\For {$i=0,\ldots, m+l$}
\State $z^{(l)}_{i+1}:=\left\{ \begin{matrix}(A-\sigma_{i}I)z^{(l)}_{i}, & i<l \\ Az^{(l)}_{i}, & i \geq l \end{matrix}\right.$ \hfill \# Compute matrix-vector product
\If {$i<l-1$}
\State $z^{(k)}_{i+1} := z^{(l)}_{i+1}, \qquad k = i+1 ,\ldots, l-1$ \hfill \text{\# Copy auxiliary basis vectors} 
\EndIf
\State \textbf{end if}
\If {$i\geq l$} \hfill \# Finalize dot-products ($g_{j,i-l+1}$)
\State $g_{j,i-l+1} := (g_{j,i-l+1}-\sum_{k=i-3l+1}^{j-1}g_{k,j}g_{k,i-l+1})/g_{j,j}; \qquad j=i-2l+2,\ldots,i-l$ \hfill \# Update transformation matrix
\State $g_{i-l+1,i-l+1}:= \sqrt{g_{i-l+1,i-l+1}-\sum_{k=i-3l+1}^{i-l}g_{k,i-l+1}^2};$
\State \# Check for breakdown and restart if required \hfill \# Square root breakdown check
\If {$i<2l$}
\State $\gamma_{i-l}:=(g_{i-l,i-l+1}+\sigma_{i-l}g_{i-l,i-l} - g_{i-l-1,i-l}\delta_{i-l-1})/g_{i-l,i-l};$ \hfill \# Add column to tridiagonal matrix
\State $\delta_{i-l}:=g_{i-l+1,i-l+1}/g_{i-l,i-l};$
\Else
\State $\gamma_{i-l}:=(g_{i-l,i-l}\gamma_{i-2l}+g_{i-l,i-l+1}\delta_{i-2l} -g_{i-l-1,i-l}\delta_{i-l-1})/g_{i-l,i-l};$
\State $\delta_{i-l}:=(g_{i-l+1,i-l+1}\delta_{i-2l})/g_{i-l,i-l};$
\EndIf
\State \textbf{end if}
\State $z^{(k)}_{i-l+k+1} := (z^{(k+1)}_{i-l+k+1} + (\sigma_k - \gamma_{i-l}) z^{(k)}_{i-l+k} - \delta_{i-l-1} z^{(k)}_{i-l+k-1})/\delta_{i-l}; \qquad k = 0, \ldots, l-1$ \hfill \# Add (auxiliary) basis vectors
\State $z^{(l)}_{i+1} := (z^{(l)}_{i+1} - \gamma_{i-l} z^{(l)}_{i}- \delta_{i-l-1} z^{(l)}_{i-1})/\delta_{i-l};$ \hfill\# Add auxiliary basis vector
\EndIf
\State \textbf{end if}
\State $g_{j,i+1}:=\left\{ \begin{matrix}(z^{(l)}_{i+1},z^{(0)}_{j}); & j=\max(0,i-2l+1),\ldots,i-l+1 \\ (z^{(l)}_{i+1},z^{(l)}_{j});  &j=i-l+2,\ldots,i+1 \end{matrix}\right.$ \hfill\# Initiate dot-products ($g_{j,i+1}$)
\If {$i=l$}
\State $\eta_{0}:=\gamma_{0};$ \quad $\zeta_{0}:={\|r_{0}\|}_2;$ \quad $p_{0}:=z^{(0)}_0/\eta_0;$
\Else \, \textbf{if} {$i\geq l+1$} \textbf{then}
\State $\lambda_{i-l}:=\delta_{i-l-1}/\eta_{i-l-1};$ \hfill \# Factorize tridiagonal matrix
\State $\eta_{i-l}:=\gamma_{i-l}-\lambda_{i-l}\delta_{i-l-1};$ 
\State $\zeta_{i-l}=-\lambda_{i-l}\zeta_{i-l-1};$ \hfill \# Compute recursive residual norm
\State $p_{i-l}=(z^{(0)}_{i-l}-\delta_{i-l-1}p_{i-l-1})/\eta_{i-l};$ \hfill \# Update search direction
\State $x_{i-l}=x_{i-l-1}+\zeta_{i-l-1}p_{i-l-1};$ \hfill \# Update approximate solution
\If {$|\zeta_{i-l}|/\|r_0\| < \tau$} RETURN; \textbf{end if}  \hfill \# Check convergence criterion
\EndIf
\EndIf
\State \textbf{end if}
\EndFor 
\end{algorithmic}
}
\end{algorithm*}

Let us expound on Alg.\,\ref{algo:PIPELCGSTAB} in some more detail.
In the $i$-th iteration of Alg.\,\ref{algo:PIPELCGSTAB} each basis $Z^{(0)}_j, Z^{(1)}_j, \ldots, Z^{(l)}_j$ is updated by adding one vector. Thus the algorithm computes a total of $l+1$ new basis vectors, i.e.: $\{v_{i-l+1} = z^{(0)}_{i-l+1}, z^{(1)}_{i-l+2}, \ldots, z^{(l-1)}_{i}, z^{(l)}_{i+1} = z_{i+1}\}$, per iteration. For each basis, the corresponding vector update in iteration $i \geq l$ is computed as follows:
\begin{align*}
\left\{ 
	\begin{array}{lll}
	z^{(0)}_{i-l+1} &\hspace{-0.2cm}= (z^{(1)}_{i-l+1} + (\sigma_{0} - \gamma_{i-l}) z^{(0)}_{i-l} - \delta_{i-l-1} z^{(0)}_{i-l-1}) / \delta_{i-l}, \\ 
	z^{(1)}_{i-l+2} &\hspace{-0.2cm}= (z^{(2)}_{i-l+2} + (\sigma_{1} - \gamma_{i-l}) z^{(1)}_{i-l+1} - \delta_{i-l-1} z^{(1)}_{i-l}) / \delta_{i-l}, \\ 
	  \quad \vdots	&\hspace{-0.2cm} \qquad \qquad\qquad \qquad \qquad	\vdots  \\
	z^{(l-1)}_{i}   &\hspace{-0.2cm}= (z^{(l)}_{i} + (\sigma_{l-1} - \gamma_{i-l}) z^{(l-1)}_{i-1} - \delta_{i-l-1} z^{(l-1)}_{i-2}) / \delta_{i-l}, \\ 
	z^{(l)}_{i+1}   &\hspace{-0.2cm}= (A z^{(l)}_{i} - \gamma_{i-l} z^{(l)}_{i} - \delta_{i-l-1} z^{(l)}_{i-1}) / \delta_{i-l}.
	\end{array}
\right.
\end{align*}
Note that all vector updates make use of the \emph{same} scalar coefficients $\gamma_{i-l}$, $\delta_{i-l}$ and $\delta_{i-l-1}$ that are computed in iteration $i$ of Alg.\,\ref{algo:PIPELCGSTAB} (lines 11-17). We also remark that only one \textsc{spmv}, namely $A z^{(l)}_{i}$, is required per iteration to compute all basis vector updates, similar to Alg.\,\ref{algo:PIPELCG}.

The merit of the intermediate basis vector updates is that they allow to replace the recurrence \eqref{eq:v_rec} for the original basis vector $v_{j+1}$ by 
relation \eqref{eq:all_z_rec}, which for $k = 0$ yields:
\begin{equation} \label{eq:all_v_rec}
v_{j+1} = z^{(0)}_{j+1} = (z^{(1)}_{j+1} + (\sigma_{0} - \gamma_{j}) v_{j} - \delta_{j-1} v_{j-1}) / \delta_{j}, 
\end{equation}
with $j \geq 0$.
Since relation \eqref{eq:all_z_relation} with $k = 0$ states that $z^{(1)}_{j+1} + \sigma_{0} v_{j} = Av_j$, expression \eqref{eq:all_v_rec} very closely resembles the finite precision Lanczos recurrence relation \eqref{eq:vbar_stab}. In particular, we point out that the matrix $G_{j+1}^{-1}$ does not occur in the recursive update \eqref{eq:all_v_rec} for $v_{j+1}$. We clarify the difference between the finite precision and exact variant of recurrence relation \eqref{eq:all_v_rec} in Section \ref{sec:rounding}.

\begin{example}
	To illustrate the methodology for constructing the basis in the stable p($l$)-CG method, consider the case where the pipeline length $l = 3$. We formulate the improved p($3$)-CG method following the derivation above. This scenario features the default Krylov subspace basis $Z_j^{(0)} = V_j$ and three auxiliary bases: $Z_j^{(1)}, Z_j^{(2)}$ and $Z_j^{(3)} = Z_j$. To compute a new vector for each auxiliary basis, the following recurrence relations are used by Alg.\,\ref{algo:PIPELCGSTAB} in iteration $i \geq 3$:
\begin{align*}
\left\{ 
	\begin{array}{lll}
	z^{(0)}_{i-2} &=& (z^{(1)}_{i-2} + (\sigma_{0} - \gamma_{i-3}) v_{i-3} - \delta_{i-4} v_{i-4}) / \delta_{i-3}, \\ 
	z^{(1)}_{i-1} 		&=& (z^{(2)}_{i-1} + (\sigma_{1} - \gamma_{i-3}) z^{(1)}_{i-2} - \delta_{i-4} z^{(1)}_{i-3}) / \delta_{i-3}, \\ 
	z^{(2)}_{i}   		&=& (z_{i} + (\sigma_{2} - \gamma_{i-3}) z^{(2)}_{i-1} - \delta_{i-4} z^{(2)}_{i-2}) / \delta_{i-3}, \\ 
	z^{(3)}_{i+1} &=& (A z_{i} - \gamma_{i-3} z_{i} - \delta_{i-4} z_{i-1}) / \delta_{i-3}.
	\end{array}
\right.
\end{align*}
The final recurrence relation to update $z^{(3)}_{i+1} = z_{i+1}$ is identical to expression \eqref{eq:z_rec} with $l= 3$. The former recurrence relation for $z^{(0)}_{i-2} = v_{i-2}$, expression \eqref{eq:v_rec}, is replaced by the (stable) relation \eqref{eq:all_v_rec}. This update explicitly uses the auxiliary variable $z^{(1)}_{i-2}$ and implicitly depends on the other $l-1 = 2$ auxiliary variables $z^{(2)}_{i-2}$ and $z^{(3)}_{i-2}$ 
through the respective recurrence relations. 
All four recurrence relations above make use of the same scalar coefficients $\delta_{i-3}, \gamma_{i-3}$ and $\delta_{i-4}$ that form the last column of the matrix $T_{i-1,i-2}$. Similar to Alg.\,\ref{algo:PIPELCG} these coefficients are computed on line 11-17 in 
Alg.\,\ref{algo:PIPELCGSTAB}, right before the recursive vector updates.
\end{example}

\begin{example}
	In the case where the pipeline length is one, i.e.~$l = 1$, the stable p($l$)-CG Alg.\,\ref{algo:PIPELCGSTAB} formally differs slightly from the original formulation of p($1$)-CG, Alg.\,\ref{algo:PIPELCG}. 
	The improved algorithm uses the following recurrence relations for $v_i$ and $z_{i+1}$ in iteration $i \geq 1$:
\begin{align*}
\left\{ 
	\begin{array}{lll}
	z^{(0)}_{i} 	&=~ (z^{(1)}_{i} + (\sigma_{0} - \gamma_{i-1}) v_{i-1} - \delta_{i-2} v_{i-2}) / \delta_{i-1}, \\ 
	z^{(1)}_{i+1} &=~ (A z_{i} - \gamma_{i-1} z_{i} - \delta_{i-2} z_{i-1}) / \delta_{i-1}.
	\end{array}
\right.
\end{align*}
where $z^{(0)}_{i} = v_{i}$ and $z^{(1)}_{i+1} = z_{i+1}$.
This implies that the \emph{only} difference between the improved and original p($1$)-CG method is the recurrence relation for $v_{i}$.
The recurrence relation for $v_{i}$ above is equivalent to the recurrence relation \eqref{eq:v_rec} in exact arithmetic, but it is numerically (more) stable as explained in Section \ref{sec:rounding}.
	
\end{example}

\subsection{Computational costs and storage requirements} \label{sec:computational}

We give an overview of implementation details of the stable p($l$)-CG method, Alg.\,\ref{algo:PIPELCGSTAB}, including global storage requirements and number of \emph{flops} (floating point operations) per iteration. We compare to the same properties for the former version of p($l$)-CG Alg.\,\ref{algo:PIPELCG} \cite{cornelis2017communication} and Ghysels' p-CG method \cite{ghysels2014hiding}. 
The latter algorithm, although mathematically equivalent to p($1$)-CG, was derived in an essentially different way. 

\begin{table}[t]
\centering
\footnotesize
\begin{tabular}{| l | c | c | c | c | c |}
\hline 
																			& \hspace{-0.15cm}\textsc{gl}\hspace{-0.1cm} &  & Flops & Time & \hspace{-0.2cm}Memory\hspace{-0.2cm} \\
																			& \hspace{-0.15cm}\textsc{sync}\hspace{-0.1cm} & \hspace{-0.15cm}\textsc{spmv}\hspace{-0.15cm} & \hspace{-0.2cm}\textsc{axpy}\&\textsc{dot}\hspace{-0.2cm} & \hspace{-0.2cm}\textsc{glred}\&\textsc{spmv}\hspace{-0.2cm} &  \\
\hline 
\hspace{-0.1cm}CG	 																	& 2 & 1 		& 10			& \hspace{-0.15cm}2 \textsc{glred} + 1 \textsc{spmv}\hspace{-0.15cm} & 3 \\
\hspace{-0.1cm}p-CG																	& 1 & 1     & 16 			& \hspace{-0.15cm}$\max$(\textsc{glred}, \textsc{spmv})\hspace{-0.15cm} & 6 \\
\hspace{-0.1cm}Alg.~\ref{algo:PIPELCG}	\hspace{-0.2cm}		& 1 & 1     & $6l+10$ & \hspace{-0.15cm}$\max$(\textsc{glred}$/l$, \textsc{spmv})\hspace{-0.15cm} & \hspace{-0.15cm}$\max$($3l+3$, $7$)\hspace{-0.15cm} \\
\hspace{-0.1cm}Alg.~\ref{algo:PIPELCGSTAB} \hspace{-0.2cm}	&	1 & 1     & $6l+10$ & \hspace{-0.15cm}$\max$(\textsc{glred}$/l$, \textsc{spmv})\hspace{-0.15cm} & \hspace{-0.15cm}$\max$($4l+1$, $7$)\hspace{-0.15cm} \\
\hline
\end{tabular}
\caption{\small Theoretical specifications of different CG variants (no preconditioning); 
`CG' denotes classic CG; `p-CG' is Ghysels' pipelined CG \cite{ghysels2014hiding}.
Columns \textsc{glsync} and \textsc{spmv} list the number of 
synchronization 
phases and \textsc{spmv}s per iteration. 
The column \emph{Flops} indicates the number of flops ($\times N$) required to compute 
\textsc{axpy}s and dot products (with $l \geq 1$). The \emph{Time} column shows the 
time spent in \textsc{glred}s (global all-reduce communications) and \textsc{spmv}s. 
\emph{Memory} counts the total number of vectors 
in memory (excl.~$x_{i-l}$ and $b$) at any time
during 
execution. 
}
\label{tab:pipelcg}
\end{table}

\subsubsection{Floating point operations per iteration} 

All Conjugate Gradient variants 
listed in Table \ref{tab:pipelcg} compute a single \textsc{spmv} in each iteration. However, as indicated by the \emph{Time} column, time per iteration may be reduced significantly by overlapping the global reduction phase with the computation of one or multiple \textsc{spmv}s. Time required by the local \textsc{axpy} and \textsc{dot-pr} computations is neglected, since these operations are assumed to scale perfectly as a function of the number of workers.

Comparing Alg.\,\ref{algo:PIPELCG} to Alg.\,\ref{algo:PIPELCGSTAB}, it is clear that the latter requires an additional $2(l-1)$ \textsc{axpy}s per iteration to update the auxiliary vectors $z^{(1)}_{i-l+2}, \ldots, z^{(l-1)}_{i}$ using the recurrence relations \eqref{eq:all_z_rec2}. However, the recurrence relation to update $v_{i-l+1}$, expression \eqref{eq:all_v_rec}, only requires $2$ \textsc{axpy} operations, instead of the $2l$ \textsc{axpy}s required to update $v_{i-l+1}$ in Alg.\,\ref{algo:PIPELCG} using expression \eqref{eq:v_rec}. Both algorithms furthermore compute $(l+1)$ local dot products (see Remark \ref{remark:dotpr}) to form the $G_j$ matrix and use two additional \textsc{axpy} operations to update the search direction $p_{i-l}$ and the iterate $x_{i-l}$. In summary, as indicated by the \emph{Flops} column in Table \ref{tab:pipelcg}, both p($l$)-CG algorithms use a total of $(6l+10)N$ flops in each iteration. The recurrence relations in the p($l$)-CG algorithm can thus be stabilized at \emph{no additional computational cost} using the framework outlined in Section \ref{sec:derivation}.

\subsubsection{Global storage requirements}

Section \ref{sec:rounding} proves that the stable p($l$)-CG method, Alg.\,\ref{algo:PIPELCGSTAB}, is extremely resilient to the presence of local rounding errors in the recurrence relations. However, this stability comes at a slightly increased storage cost compared to Alg.\,\ref{algo:PIPELCG}. The latter requires to store 
$2l+1$ vectors of the $V_j$ basis (required for vector updates), $l+1$ vectors of the auxiliary basis $Z_j$ (for vector updates and dot product computations, see also Remark \ref{remark:dotpr}), and the vector $p_{i-l}$ at any time during the execution of the algorithm 
from iteration $i \geq 2l-1$ onward. 
In contrast, Alg.\,\ref{algo:PIPELCGSTAB} stores the three most recently updated vectors in each of the bases $Z^{(0)}_j, \ldots, Z^{(l)}_j$ (which include the bases $V_j$ and $Z_j$). In addition, $l$ vectors in the $Z_j$ basis need to be stored for dot product computations. 
Table \ref{tab:pipelcg} summarizes the storage requirements for different variants of the CG algorithm. Alg.\,\ref{algo:PIPELCG} keeps a total of $3l+2$ vectors in memory at any time during execution, whereas Alg.\,\ref{algo:PIPELCGSTAB} stores $4l+1$ vectors. 
The memory overhead for the stable p($l$)-CG method thus amounts to a modest 
$l-2$ vectors.

\subsection{A stable preconditioned p($l$)-CG algorithm} \label{sec:preconditioning}

Preconditioning is 
indispensable to efficiently solve linear systems in practice. We 
briefly comment on the straightforward extension of Alg.\,\ref{algo:PIPELCGSTAB} to include a preconditioner. 
This section follows the standard methodology that was described for 
pipelined CG in \cite{ghysels2014hiding} and for Alg.\,\ref{algo:PIPELCG} in 
\cite{cornelis2017communication}. 

Let the preconditioner be given by the matrix $M^{-1}$. We aim to solve the left-preconditioned linear system $M^{-1}Ax=M^{-1}b$, where $M$ and $A$ are both symmetric positive definite matrices. This assumption does not necessarily imply that $M^{-1}A$ is symmetric. Nonetheless, symmetry can be preserved by observing that $M^{-1}A$ is self-adjoint with respect to the $M$ inner product $(x,y)_{M} = (Mx,y) = (x,My)$. The basic strategy is thus to replace all Euclidean dot products occurring in Alg.\,\ref{algo:PIPELCGSTAB} with $M$ inner products. 

We cannot simply use the matrix $M$ to calculate these $M$ inner products, since the preconditioner inverse is in general not available. However, by introducing the \textit{unpreconditioned auxiliary variables} $u_{i}=M z^{(l)}_{i}$ and observing that these variables again satisfy a Lanczos type 
relation: 
\begin{equation} \label{eq:zhatid3}
u_{i+1} 
= \left\{ \begin{matrix}Az^{(l)}_{i} -\sigma_{i}u_{i}  & i<l, \\ (Az^{(l)}_{i} - \gamma_{i-l} u_{i}- \delta_{i-l-1} u_{i-1})/\delta_{i-l} & i \geq l,  \end{matrix}\right.
\end{equation}
this obstacle is circumvented. Using these unpreconditioned auxiliary variables $u_{i+1}$, the dot products $g_{j,i+1}$ for $0 \leq j \leq i+1$ can be computed as follows. 
For $i < l$ it holds that
\begin{equation} \label{eq:gdotpr1prec}
g_{j,i+1} = (z^{(l)}_{i+1},z^{(l)}_{j})_M = (M z^{(l)}_{i+1},z^{(l)}_{j}) = (u_{i+1},z^{(l)}_{j}) , 
\end{equation}
with $j=0,\ldots,i+1$.
For $i \geq l$ we find
\begin{equation} \label{eq:gdotpr2prec}
g_{j,i+1}=\left\{ 
  \begin{matrix}
	  (z^{(l)}_{i+1},z^{(0)}_{j})_M = (u_{i+1},z^{(0)}_{j}), \qquad \qquad \qquad \qquad \qquad \qquad \vspace{0.1cm}\\
		\quad \text{for}~ j=\max(0,i-2l+1),\ldots,i-l+1, \vspace{0.1cm} \\
	  (z^{(l)}_{i+1},z^{(l)}_{j})_M = (u_{i+1},z^{(l)}_{j}), \qquad \qquad \qquad \qquad \qquad \qquad \vspace{0.1cm}\\
		\quad \text{for}~ j=i-l+2,\ldots,i+1. \qquad \qquad \qquad \qquad 
	\end{matrix} 
\right. 
\end{equation}
This allows to formulate a preconditioned version of Alg.\,\ref{algo:PIPELCGSTAB} by adding the recurrence relation \eqref{eq:zhatid3} for the unpreconditioned auxiliary variables $u_{i}$, and replacing the dot products on line $23$ and $25$ by expressions \eqref{eq:gdotpr1prec} and \eqref{eq:gdotpr2prec} respectively. From an implementation point of view the extension to the preconditioned algorithm only requires the application of the preconditioner $M^{-1}$, two additional \textsc{axpy} operations and 
storage of three additional vectors 
in memory.

\subsection{Rounding error analysis for the improved method} \label{sec:rounding} 

We now consider the finite precision equivalents of the above recurrence relations for all basis vectors in the improved version of p($l$)-CG, Alg.\,\ref{algo:PIPELCGSTAB}. The actually computed finite precision variants of the basis vectors and scalar variables are again denoted by 
bars.
 
Let the true basis vector $\bar{\bold{v}}_j$ satisfy the actual Lanczos relation without local roundoff (for $0 \leq j < i-l$)
\begin{equation} \label{eq:rec_v_bold}
\bar{\bold{v}}_{j+1} = ( A\bar{v}_j - \bar{\gamma}_j \bar{v}_j - \bar{\delta}_{j-1} \bar{v}_{j-1} ) / \bar{\delta}_j,  
\end{equation}
and let, analogously, the true auxiliary basis vectors be defined (for $j \geq k$ and $0 \leq k \leq l$) as 
\begin{equation} \label{eq:all_z_lanczos}
\bar{\bold{z}}^{(k)}_{j+1} = (A \bar{z}^{(k)}_{j} - \bar{\gamma}_{j-k} \bar{z}^{(k)}_{j} - \bar{\delta}_{j-k-1} \bar{z}^{(k)}_{j-1}) / \bar{\delta}_{j-k}.
\end{equation}
On the other hand, the finite precision variants of the recurrence relations \eqref{eq:all_z_rec2}, which are actually computed in Alg.\,\ref{algo:PIPELCGSTAB}, are for $j \geq k$ and $0 \leq k < l$ given by the expressions
\begin{equation} \label{eq:all_z_fin}
\bar{z}^{(k)}_{j+1} = (\bar{z}^{(k+1)}_{j+1} + (\sigma_{k} - \bar{\gamma}_{j-k}) \bar{z}^{(k)}_{j} - \bar{\delta}_{j-k-1} \bar{z}^{(k)}_{j-1}) / \bar{\delta}_{j-k} + \xi^{(k)}_{j+1},
\end{equation}
whereas for $k = l$ the following finite precision relation holds (with $j \geq l$):
\begin{equation} \label{eq:all_z_last}
\bar{z}^{(l)}_{j+1} = (A\bar{z}^{(l)}_{j} - \bar{\gamma}_{j-l} \bar{z}^{(l)}_{j} - \bar{\delta}_{j-l-1} \bar{z}^{(l)}_{j-1}) / \bar{\delta}_{j-l} + \xi^{(l)}_{j+1}, 
\end{equation}
The local rounding errors $\xi^{(k)}_{j+1}$ $(0 \leq k \leq l)$ in the above expressions can be bounded by machine precision $\epsilon$ multiplied by the norms of the respective vectors, i.e.~$\|\xi^{(k)}_{j+1}\| \leq 
\mathcal{O}(\epsilon)$.

\begin{figure*}
    \centering
    \begin{subfigure}[b]{0.49\textwidth}
				\includegraphics[width=\textwidth]{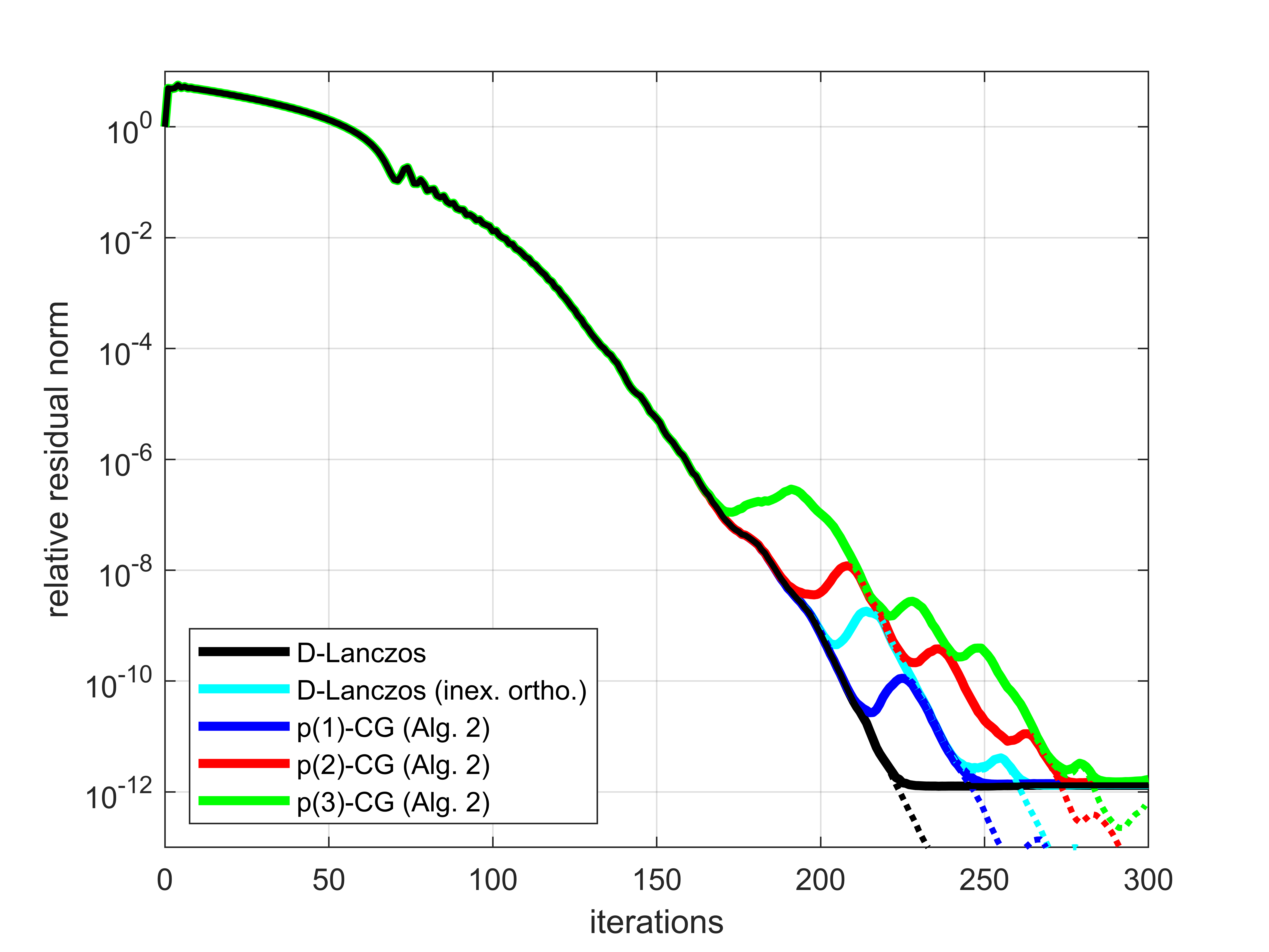} 
        \caption{}
        \label{fig:figure2a}
    \end{subfigure}
    \begin{subfigure}[b]{0.49\textwidth}
				\includegraphics[width=\textwidth]{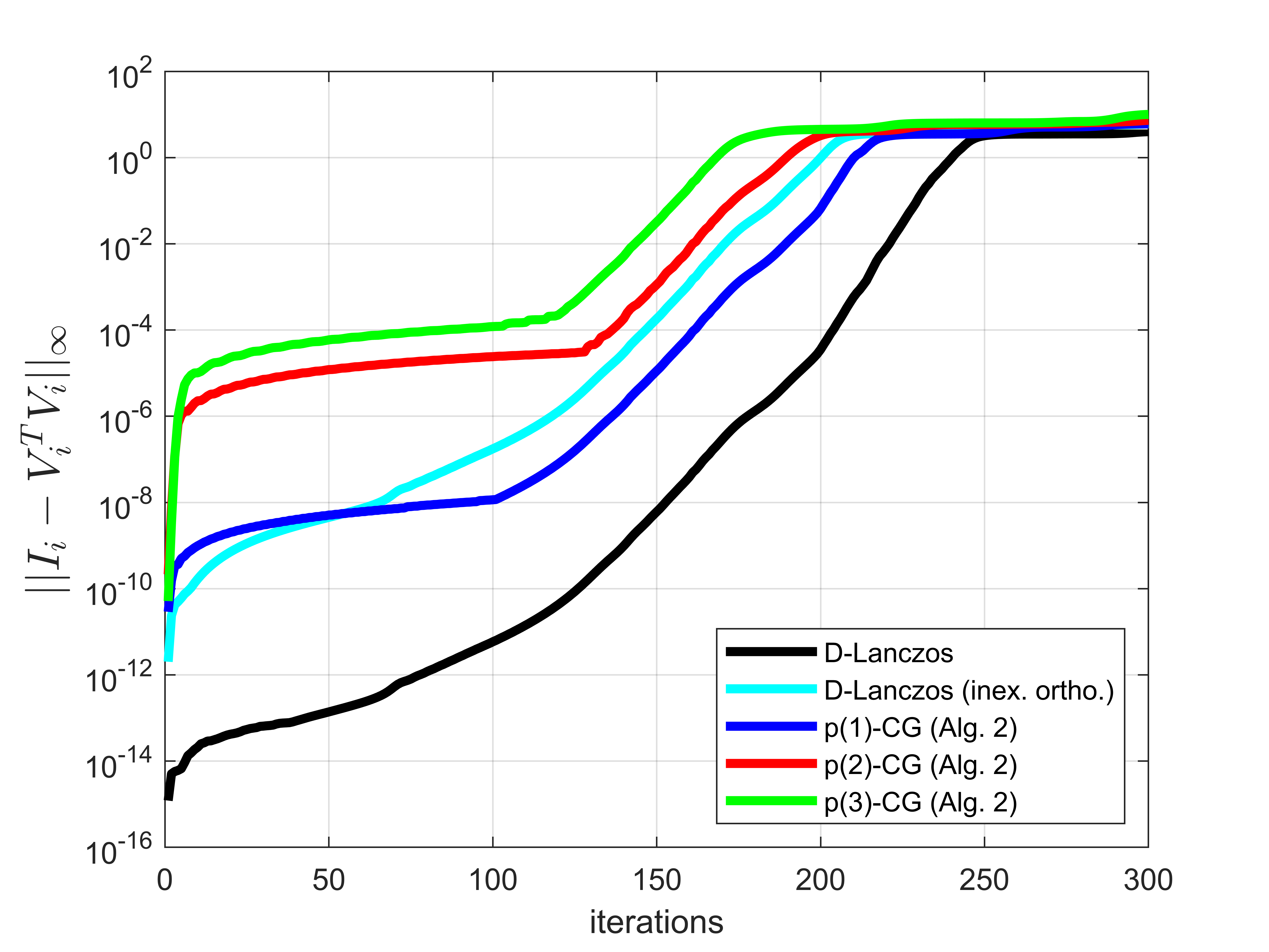}
        \caption{}
        \label{fig:figure2b}
    \end{subfigure}
		\vspace{-0.2cm}
\caption{\textbf{(B1)} \small
Fig.\,\ref{fig:figure1} analogue. (a) 
True residual norm $\|b-A\bar{x}_i\|/\|b\|$ (full line) and recursive residual norm $\|\bar{r}_i\|/\|b\|$ (dotted line) for different CG variants for a 2D $5$-point stencil discretized $100 \times 100$ unknowns Laplace problem.  
Recurrence relations \eqref{eq:all_z_rec2}-\eqref{eq:specific_z_rec} are used in stable p($l$)-CG method, Alg.\,\ref{algo:PIPELCGSTAB}. 
(b) Basis orthogonality characterized by $\|I - \bar{V}_i^T \bar{V}_i\|_{\infty}$. 
}
\label{fig:figure2}
\end{figure*}

Switching to matrix forms to simplify the notation, we write the Lanczos relations \eqref{eq:all_z_lanczos} as
\begin{equation} \label{eq:all_z_lanczos_mat}
	\bar{\bold{Z}}^{(k)}_{j+1} - \bar{Z}^{(k)}_{j+1} = (A\bar{Z}^{(k)}_{j} - \bar{Z}^{(k)}_{j+1} \bar{T}_{j+1,j}) \bar{\Delta}^+_{j+1,j}, 
\end{equation}
for $j > k$ and $0 \leq k \leq l$.
Note that this expression neglects the ``\emph{initial}'' gaps, i.e.~the quantities $\bar{\bold{z}}^{(k)}_{j} - \bar{z}^{(k)}_{j}$ for $0 \leq j \leq k$, consisting of the local rounding errors from the application of the matrix polynomials $P_j(A)$ to the initial vector $v_0$, which are computed explicitly in the stable p($l$)-CG method, Alg.\,\ref{algo:PIPELCGSTAB} (line 3).
The recurrence relations \eqref{eq:all_z_fin} are summarized by the matrix expressions
\begin{equation} \label{eq:all_z_fin_mat}
	\bar{Z}^{(k+1)}_{2:j+1} = \bar{Z}^{(k)}_{j+1} \bar{T}_{j+1,j} - \sigma_k \bar{Z}^{(k)}_{j} - \Xi^{(k)}_{j+1} \bar{\Delta}_{j+1,j},  
\end{equation}
with $j > k$ and $0 \leq k < l$,
where $\Xi^{(k)}_{j+1} = [\xi^{(k)}_0,\ldots,\xi^{(k)}_{j}]$ and $\bar{Z}^{(k+1)}_{2:j+1} = [\bar{z}^{(k+1)}_{1}, \bar{z}^{(k+1)}_{2},\ldots, \bar{z}^{(k+1)}_{j}]$.
The recurrence relation \eqref{eq:all_z_last} for $\bar{z}^{(l)}_{j+1}$ can be written as
\begin{equation} \label{eq:all_z_last_mat}
	A\bar{Z}^{(l)}_{j} = \bar{Z}^{(l)}_{j+1} \bar{T}_{j+1,j} - \Xi^{(l)}_{j+1} \bar{\Delta}_{j+1,j},
\end{equation}
with $j > l$.
To find an expression for the gap on the basis vectors in $\bar{V}_j$, i.e.~the gap $\bar{\bold{Z}}^{(0)}_{j} - \bar{Z}^{(0)}_{j}$, we now progressively compute the gaps on the auxiliary bases. Starting from $\bar{\bold{Z}}^{(l)}_{j} - \bar{Z}^{(l)}_{j}$, we compute $\bar{\bold{Z}}^{(l-1)}_{j} - \bar{Z}^{(l-1)}_{j}$, in which we substitute the gap on $\bar{Z}^{(l)}_{j}$, followed by $\bar{\bold{Z}}^{(l-2)}_{j} - \bar{Z}^{(l-2)}_{j}$, etc., until we eventually obtain an expression for the gap $\bar{\bold{Z}}^{(0)}_{j} - \bar{Z}^{(0)}_{j}$.
For $k = l$ we derive from \eqref{eq:all_z_lanczos_mat} and \eqref{eq:all_z_last_mat} that
\begin{equation} \label{eq:gap_particular_l}
	\bar{\bold{Z}}^{(l)}_{j+1} - \bar{Z}^{(l)}_{j+1} = - \Xi^{(l)}_{j+1},
	\quad j > l.
\end{equation}
As indicated by this expression, the gaps on the basis vectors in the basis $\bar{Z}^{(l)}_{j}$ thus consist of local rounding errors only. The following auxiliary lemma can easily be proven by induction.
\begin{lemma} \label{lemma:lemma1}
Let $\bar{Z}^{(k)}_j$ be defined by \eqref{eq:all_z_fin_mat} and $\bar{\bold{Z}}^{(k+1)}_{j+1}$ be defined by \eqref{eq:all_z_lanczos_mat}. Then it holds
for  $ j > k$ and $0 \leq k < l$ that
\begin{equation}
	A \bar{Z}^{(k)}_j - \sigma_k \bar{Z}^{(k)}_j = \bar{\bold{Z}}^{(k+1)}_{2:j+1} + \Theta^{(k+1)}_{2:j+1},
\end{equation}
where $\Theta^{(k+1)}_{2:j+1} = [\theta^{(k+1)}_{1} , \theta^{(k+1)}_{2}, \ldots , \theta^{(k+1)}_{j}]$ with $\theta^{(k+1)}_{j}$ a local rounding error that is bounded by $\mathcal{O}(\epsilon)$ and $\bar{\bold{Z}}^{(k+1)}_{2:j+1} = [\bar{\bold{z}}^{(k+1)}_{1}, \bar{\bold{z}}^{(k+1)}_{2},\ldots, \bar{\bold{z}}^{(k+1)}_{j}]$.
\end{lemma}
\noindent Next, we combine expressions \eqref{eq:all_z_lanczos_mat} and \eqref{eq:all_z_fin_mat} and Lemma \ref{lemma:lemma1} for the case $k = l-1$. We obtain (with $j > l-1$)
\begin{align*} 
	&\bar{\bold{Z}}^{(l-1)}_{j+1} - \bar{Z}^{(l-1)}_{j+1} \\
	&~~= (A\bar{Z}^{(l-1)}_{j} - \bar{Z}^{(l)}_{2:j+1} - \sigma_{l-1} \bar{Z}^{(l-1)}_{j} ) \bar{\Delta}^+_{j+1,j} - \Xi^{(l-1)}_{j+1}   \\
	&~~= (\bar{\bold{Z}}^{(l)}_{2:j+1} - \bar{Z}^{(l)}_{2:j+1}  + \Theta^{(l)}_{2:j+1}) \bar{\Delta}^+_{j+1,j}  - \Xi^{(l-1)}_{j+1}.
\end{align*}
Hence, the gaps on the basis vectors $\bar{Z}^{(l-1)}_{j}$ are coupled to the gaps on the basis vectors $\bar{Z}^{(l)}_{j}$. After substitution of expression \eqref{eq:gap_particular_l} it is clear that $\bar{\bold{Z}}^{(l-1)}_{j+1} - \bar{Z}^{(l-1)}_{j+1}$ consists only of local rounding errors. 
The above relation can be generalized for any $k \in \{0, 1, \ldots, l-1\}$ as follows:
\begin{align} \label{eq:gap_general_k}
	&\bar{\bold{Z}}^{(k)}_{j+1} - \bar{Z}^{(k)}_{j+1} \notag \\
	&~~= (A\bar{Z}^{(k)}_{j} - \bar{Z}^{(k+1)}_{2:j+1} - \sigma_{k} \bar{Z}^{(k)}_{j}  ) \bar{\Delta}^+_{j+1,j} - \Xi^{(k)}_{j+1} \notag \\
	&~~= (\bar{\bold{Z}}^{(k+1)}_{2:j+1} - \bar{Z}^{(k+1)}_{2:j+1} + \Theta^{(k+1)}_{2:j+1}) \bar{\Delta}^+_{j+1,j} - \Xi^{(k)}_{2:j+1}, 
\end{align}
where $j > k$.
After subsequent substitution of this expression starting from $\bar{\bold{Z}}^{(0)}_{j+1} - \bar{Z}^{(0)}_{j+1}$ up until the characterization of $\bar{\bold{Z}}^{(l)}_{j+1} - \bar{Z}^{(l)}_{j+1}$, see \eqref{eq:gap_particular_l}, it appears that the gap $\bar{\bold{Z}}^{(0)}_{j+1} - \bar{Z}^{(0)}_{j+1} = \bar{\bold{V}}_{j+1} - \bar{V}_{j+1}$ on the original Krylov subspace basis 
is just a sum of local rounding errors. 
No rounding error propagation takes place in the stable p($l$)-CG method, Alg.\,\ref{algo:PIPELCGSTAB}, on any of the (auxiliary) basis vectors, 
see expressions \eqref{eq:gap_particular_l}-\eqref{eq:gap_general_k}. 
By introducing the intermediate auxiliary bases $\bar{Z}^{(1)}_j, \ldots, \bar{Z}^{(l-1)}_j$ for the recursive computation of the basis $\bar{V}_j$, the dependency of the basis vector gaps on the possibly ill-conditioned matrix $\bar{G}^{-1}_j$, cf.~expressions \eqref{eq:Z_BAR} and \eqref{eq:gap_plcg2}, is thus removed, resulting in a numerically stable algorithm.

\section{Numerical results} \label{sec:numerical} 

We present various numerical experiments to benchmark the stable p($l$)-CG method 
proposed in Section \ref{sec:stabilizing}. The benchmark problems exemplify both the 
performance of the improved p($l$)-CG method on large scale parallel hardware 
as well as its 
error resilience 
compared to other 
CG variants.
Performance measurements 
result from a 
PETSc \cite{petsc-web-page} implementation of the 
p($l$)-CG algorithm on a distributed memory machine using the message passing
paradigm (MPI). 

\subsection{Hardware and software specifications} \label{sec:hardware}

Parallel performance experiments 
are performed on up to 128 compute nodes
of a cluster consisting of two 14-core Intel E5-2680v4 Broadwell generation CPUs 
each (28 cores per node). Nodes are connected through an EDR InfiniBand network. 
We use PETSc version 3.8.3 \cite{petsc-web-page}. The MPI library used for this 
experiment is Intel MPI 2018v3. The PETSc environment variables 
\texttt{MPICH\_ASYNC\_PROGRESS=1} and
\texttt{MPICH\_MAX\_THREAD\_SAFETY=multiple} 
ensure optimal parallelism by allowing 
for asynchronous non-blocking global communication. 
Timing results reported in this manuscript are the most favorable results 
(
smallest overall run-time) over 3 individual runs of each method.
Experiments also show results for Ghysels' p-CG method \cite{ghysels2014hiding} as a reference. The p-CG method is similar to 
p($1$)-CG in operational cost (see Table \ref{tab:pipelcg}), but features a significant 
loss of attainable accuracy due to rounding error propagation in its recurrence relations \cite{cools2018analyzing}, 
similar to p($l$)-CG, Alg.\,\ref{algo:PIPELCG}.

\subsection{Benchmark (B1): 2D Laplace PDE}

Fig.\,\ref{fig:figure2} is the analogue of the experiment reported in Fig.\,\ref{fig:figure1} where the p($l$)-CG Alg.\,\ref{algo:PIPELCG} is replaced by the improved Alg.\,\ref{algo:PIPELCGSTAB}. The model problem solved is the 2D Laplace equation
\begin{equation} \label{eq:laplace_eq}
	- \Delta u = b, \qquad 0 \leq x,y \leq 1,
\end{equation}
with homogeneous Dirichlet boundary conditions, discretized using second order finite differences on a uniformly spaced $100 \times 100$ grid.
In contrast to Fig.\,\ref{fig:figure1a}, no loss of accuracy due to local rounding error amplification in the recurrence relations is observed in Fig.\,\ref{fig:figure2a}. The stabilized recurrences \eqref{eq:all_z_rec2}-\eqref{eq:specific_z_rec} ensure that the quantity $\|(A\bar{V}_i - \bar{V}_{i+1}\bar{T}_{i+1,i}) \bar{U}_i^{-1} \bar{q}_i\|$ is of order machine precision. The true residual norms (full lines) and recursively computed residual norms (dotted lines) coincide up to $\|b-A\bar{x}_i\|_2/\|b\|_2 = 1.0$e-$12$ (maximal attainable accuracy) for all methods in Fig.\,\ref{fig:figure2a}. Fig.\,\ref{fig:figure2b} quantifies the inexact orthogonality for the different methods 
which is comparable to Fig.\,\ref{fig:figure1b}. Hence, similarly to Fig.\,\ref{fig:figure1a} a delay of convergence 
may be observed in Fig.\,\ref{fig:figure2a}.


\begin{figure*}
    \centering
    \begin{subfigure}[b]{0.39\textwidth}
				\includegraphics[width=\textwidth]{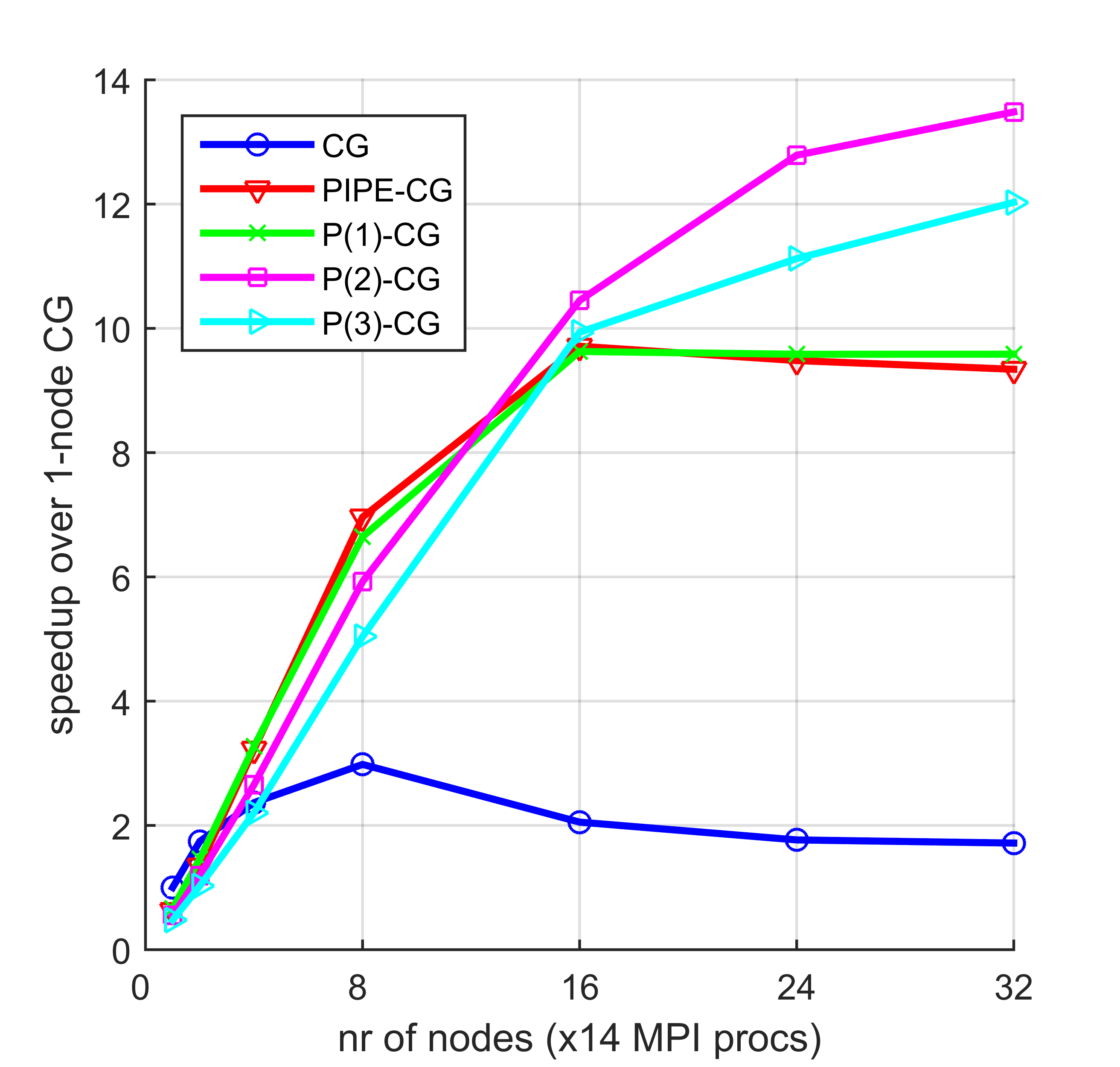} 
        \caption{}
        \label{fig:figure3a}
    \end{subfigure}
		\qquad \qquad
    \begin{subfigure}[b]{0.39\textwidth}
				\includegraphics[width=\textwidth]{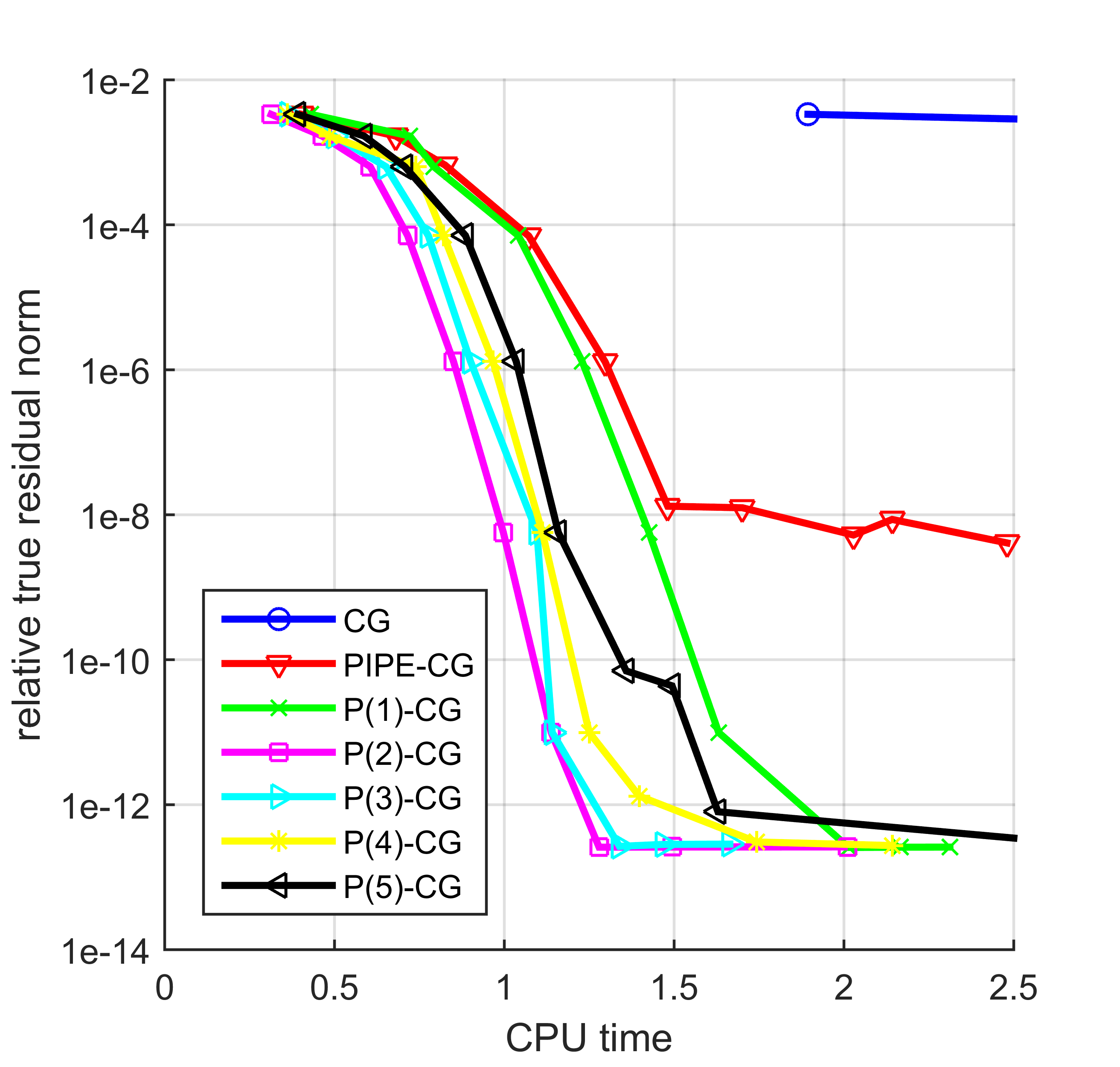} 
        \caption{}
        \label{fig:figure3b}
    \end{subfigure}
		\vspace{-0.3cm}
\caption{\textbf{(B1)} \small
(a) Strong scaling experiment on up to $32$ nodes ($448$ processes), 2D $5$-point stencil 
Laplace problem (PETSc KSP ex2), $1,750\times1,750$ unknowns. No preconditioner. 
All methods converged to $\|b-A\bar{x}_i\|_2/\|b\|_2 = 6.3$e-$4$ 
in 
$1,500$ iterations.
(b) Accuracy experiment on $32$ nodes. 
Relative residual norm $\|b-A\bar{x}_i\|_2/\|b\|_2$ 
as a function of total time spent (in s.). 
\label{fig:figure3}
} 
\end{figure*}

Fig.\,\ref{fig:figure3} shows a performance experiment on the hardware and software setup specified above.
A linear system resulting from discretization of the 2D Laplace equation \eqref{eq:laplace_eq} with exact solution $\hat{x}_j = 1$, right-hand side $b = A\hat{x}$ 
and initial guess $x_0 = 0$ is solved. This problem is available as example $2$ in the PETSc 
Krylov subspace solvers (KSP) folder. The simulation domain is discretized using a $1,750\times1,750$ 
uniform finite difference mesh ($3,062,500$ unknowns). No preconditioner is applied. 
For p($l$)-CG Chebyshev shifts $\{ \sigma_0, \ldots, \sigma_{l-1} \}$ are used based on 
the interval $[\lambda_{\min}, \lambda_{\max}] = [0,8]$ (known analytically), see \eqref{eq:chebyshev}.

Fig.\,\ref{fig:figure3a} shows the speedup achieved by different CG methods over single-node classic CG for various pipeline lengths and node setups.
Classic CG scales poorly for this model problem; no speedup is achieved beyond $8$ nodes. The pipelined methods scale well.
The length one p-CG and p($1$)-CG method achieve a relative speedup of approximately $5\times$ compared to classic CG when both 
are executed on $16$ nodes. The longer pipelined p($2$)-CG and p($3$)-CG methods out-scale the latter method, with p($2$)-CG 
achieving a $7\times$ speedup relative to classic CG on $32$ nodes. When $l = 1$ the global communication phase in each iteration is 
only partially `hidden' behind the \textsc{spmv} computation, whereas overlapping with more than two \textsc{spmv}s
by using pipelines length $l \geq 3$ does not seem to improve performance further. Pipeline length $l=2$ 
is optimal for this problem,
striking a good balance between overlapping communication and introducing additional auxiliary vectors.

Fig.\,\ref{fig:figure3b} 
plots the relative residual norms as a function of the total time spent (in s.) by various CG algorithms (on 500 iteration intervals). 
The p-CG method stagnates around $\|b-A\bar{x}_i\|_2/\|b\|_2 = 1.0$e-$8$ and is unable to attain a better accuracy regardless of computational 
effort. 
The stable p($l$)-CG methods all are able to attain a much higher accuracy,
stagnating around $\|b-A\bar{x}_i\|_2/\|b\|_2 = 2.7$e-$13$ for $l \in \{1,2,3,4,5\}$. Maximal attainable accuracy is reached in 2.01 s.~for p($1$)-CG,
1.28 s.~for p($2$)-CG, 1.34 s.~for p($3$)-CG, 1.74 s.~for p($4$)-CG, 2.73 s.~for p($5$)-CG. Note that classic CG 
attains $\|b-A\bar{x}_i\|_2/\|b\|_2 = 2.5$e-$13$ in $13.9$ s.~(outside the graph). These results are consistent with Fig.\,\ref{fig:figure3a}, 
indicating that a pipeline length of $l = 2$ suffices to hide the global communication phase for this problem.


\subsection{Benchmark (B2): 3D Hydrostatic Ice Sheet Flow}


\begin{figure*}
    \centering
    \begin{subfigure}[b]{0.39\textwidth}
				\includegraphics[width=\textwidth]{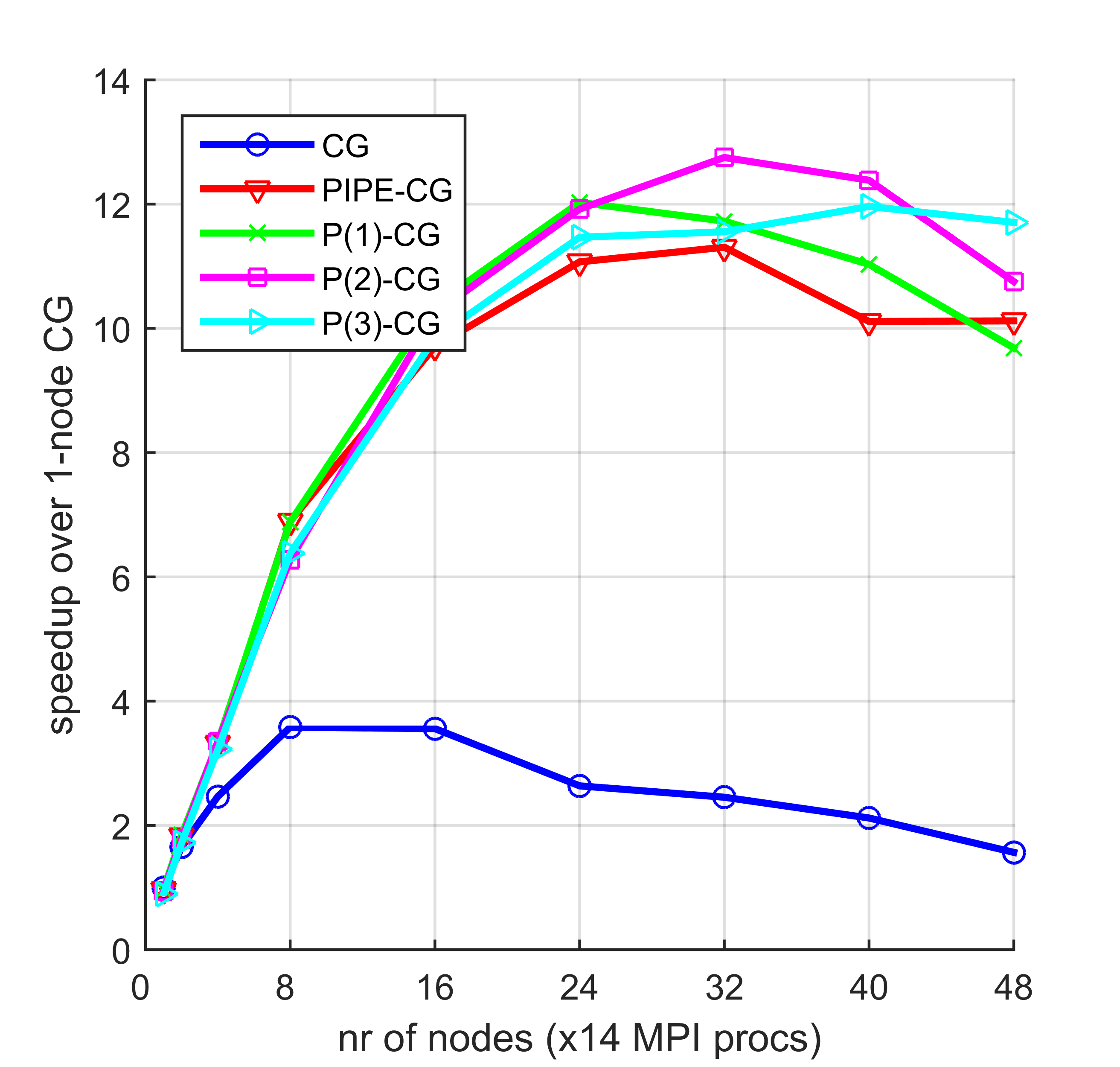} 
        \caption{}
        \label{fig:figure4a}
		    \end{subfigure}
				\qquad \qquad
    \begin{subfigure}[b]{0.39\textwidth}
				\includegraphics[width=\textwidth]{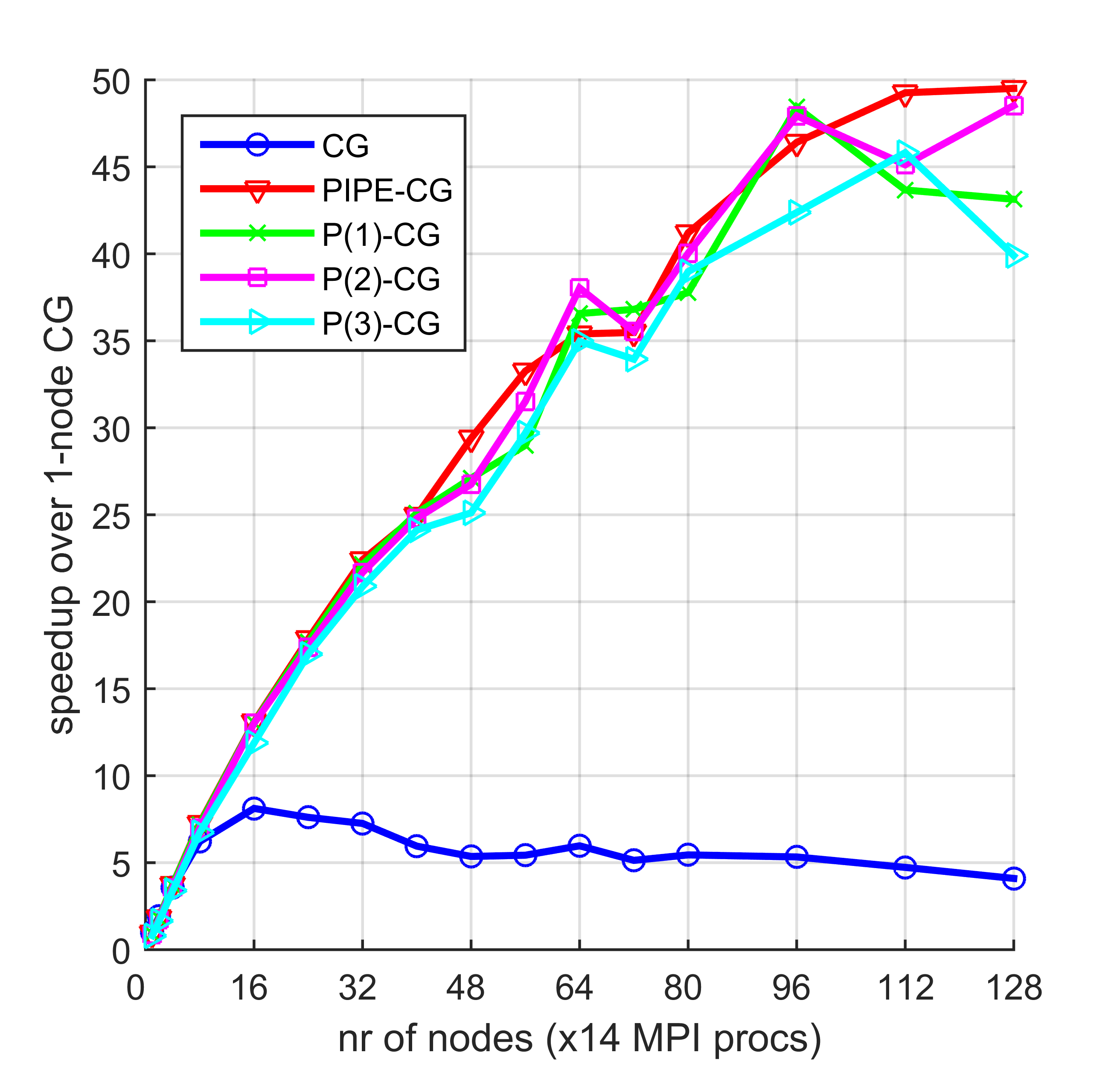} 
        \caption{}
        \label{fig:figure4b}
    \end{subfigure}
		\begin{subfigure}[b]{0.98\textwidth}
				\includegraphics[width=\textwidth]{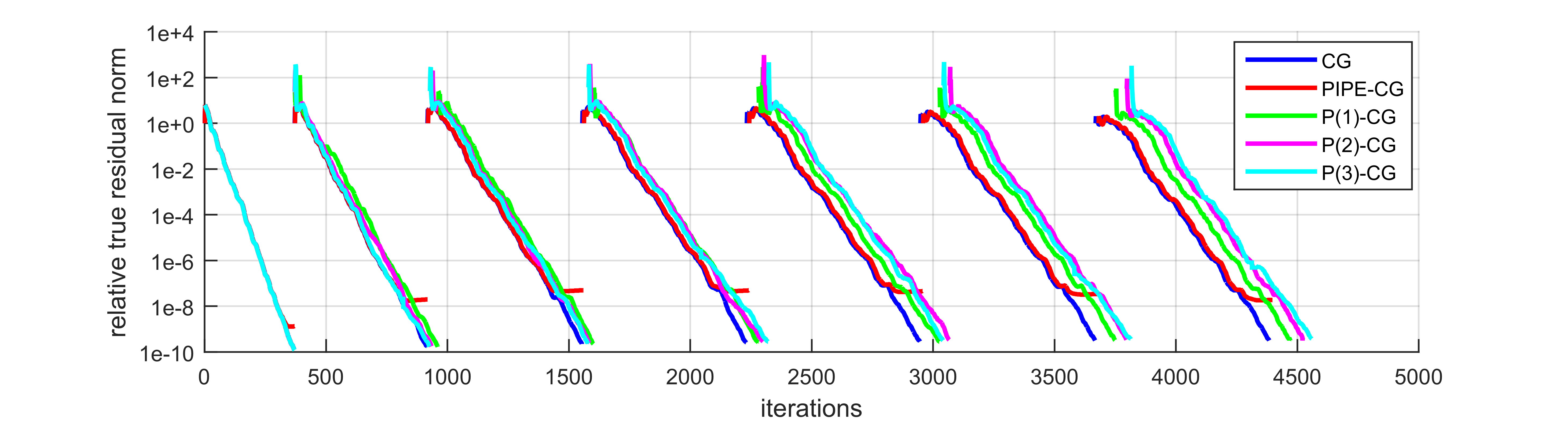} 
        \caption{}
        \label{fig:figure4c}
    \end{subfigure}
		\vspace{-0.2cm}
\caption{\textbf{(B2)} \small
(a) Non-linear 3D hydrostatic (Blatter/Pattyn) equations for ice sheet flow (PETSc SNES ex48), 
$100\times100\times20$ finite elements. Block Jacobi preconditioner. 
Inner (Krylov) solver tolerance: 
$\|b-A\bar{x}_i\|_2/\|b\|_2 \leq 1.0$e-$10$. 
(b) Non-linear 3D hydrostatic equations for ice sheet flow,
$150\times150\times100$ finite elements. 
Inner solver tolerance: 
$\|b-A\bar{x}_i\|_2/\|b\|_2 \leq 1.0$e-$5$.
(c) Accuracy experiment on $8$ nodes ($112$ processes) for the setup used in (b). 
Residual norm $\|b-A\bar{x}_i\|_2/\|b\|_2$ 
as a function of iterations 
(inner 
tolerance = $1.0$e-$10$, outer 
tolerance = $1.0$e-$8$).}
\label{fig:figure4}
\end{figure*}

Fig.\,\ref{fig:figure4a} shows the result of two strong scaling experiments for the 3D Hydrostatic Ice Sheet Flow problem, see \cite{brown2013achieving} for a full problem description. 
An implementation of this problem is available as example $48$ is the PETSc Scalable Nonlinear Equations Solvers (SNES) folder. The Blatter/Pattyn equations 
are discretized using $100\times100\times20$ (Fig.\,\ref{fig:figure4a}) and $150\times150\times100$ (Fig.\,\ref{fig:figure4b}) finite elements respectively. 
Hard- and software specifications are presented in Section \ref{sec:hardware}. A Newton-Krylov outer-inner iteration is used to solve the non-linear problem. The CG methods used as inner solver are combined with a block Jacobi preconditioner\footnote{\textbf{Note:} To attain good parallel scaling using pipelined Krylov 
methods it is generally beneficial to choose a preconditioner which does not require global communication. We hence use (block) stationary iterative methods as preconditioners that only use \textsc{spmv}s and \textsc{axpy}s but avoid computing dot products in both examples \textbf{(B2)} and \textbf{(B3)}.} (one block Jacobi step per CG iteration; one block per processor). The relative tolerance of the inner solvers is set to $\|b-A\bar{x}_i\|_2/\|b\|_2 = 1.0$e-$10$ (Fig.\,\ref{fig:figure4a}) and $\|b-A\bar{x}_i\|_2/\|b\|_2 = 1.0$e-$5$ (Fig.\,\ref{fig:figure4b}), while the outer relative tolerance is chosen to be $1.0$e-$8$. Seven Newton iterations are needed to reach the outer tolerance. Chebyshev shifts are based on the interval $[\lambda_{\min}, \lambda_{\max}] = [0,2]$, which is chosen in an \emph{ad hoc} fashion for this problem based on the presumed clustering of the spectrum around $1$. 

In Fig.\,\ref{fig:figure4b} the scalability of p($l$)-CG for $l \in \{1,2,3\}$ is comparable to that of the p-CG method. On $128$ nodes a speedup of approximately $10 \times$ is measured for the pipelined methods compared to classic CG on $128$ nodes. There is no gain in using longer pipeline lengths for this problem and hardware setup, indicating that the amount of computational work of a single iteration (\textsc{spmv} + \textsc{prec}) suffices to hide the communication in the global reduction in each iteration. It is 
expected 
that longer pipeline lengths 
out-scale the methods with shorter pipelines on very large numbers of nodes, since communication costs would increase accordingly. This is illustrated to some extend by the smaller 
problem reported in Fig.\,\ref{fig:figure4a}, which shows that for heavily communication bound problems the use of longer pipelines is beneficial for this benchmark. 
Fig.\,\ref{fig:figure4a} indicates that depending on the amount of hardware parallelism (
number of nodes) p($1$)-CG (4-24 nodes), p($2$)-CG (32-40 nodes) or p($3$)-CG (48+ nodes) respectively display the biggest speedup. 

An accuracy experiment for the $150\times150\times100$ FE benchmark problem is shown in Fig.\,\ref{fig:figure4c}. The experiment is run on $8$ nodes and the relative tolerance of the inner solver is 
$\|b-A\bar{x}_i\|_2/\|b\|_2 = 1.0$e-$10$. A relatively small delay of convergence is observed for the p($l$)-CG method, with the effect worsening for longer pipelines.
The convergence delay is due to restarts (caused by a square root breakdown) which come forth from the ill-conditioned auxiliary basis $\bar{Z}_{j}$. The effect is negligible since the maximal total delay after seven Newton iterations is $174$ iterations (comparing p($3$)-CG to CG), relative to a total of 
$\sim 4,500$ Krylov iterations. 
The p-CG method fails to reach the inner tolerance $1.0$e-$10$. 
As shown by the analysis in Section \ref{sec:rounding}, no rounding error propagation occurs in the stable p($l$)-CG algorithm and the method 
is able to attain a highly accurate solution. 


\subsection{Benchmark (B3): 2D Bratu Solid Fuel Ignition}

The 2D Bratu Solid Fuel Ignition problem results from a finite difference discretization of the nonlinear PDE 
\begin{equation}
	\Delta u - \lambda_{B} \, e^u = 0, \qquad 0 < x, y < 1,
\end{equation}
with homogeneous Dirichlet boundary conditions on a uniformly spaced 2D grid.
Its implementation can be found in the PETSc SNES folder as example 5. The two-dimensional domain is discretized with $1,000 \times 1,000$ uniformly spaced grid points. The Bratu parameter $\lambda_{B}$ is set to $6.0$, implying significant non-linearity. Different CG methods are used as the inner solver in the Newton-Krylov scheme. The preconditioner is a SOR-preconditioned Chebyshev iteration (three SOR steps per Chebyshev iteration; one Chebyshev step per CG iteration).  Outer solver relative tolerance is set to $1.0$e-$8$, which is attained by all methods after three Newton steps. For p($l$)-CG Chebyshev shifts based on the interval $[\lambda_{\min}, \lambda_{\max}] = [0, 1.2]$ are used. A relative tolerance of $\|b-A\bar{x}_i\|_2/\|b\|_2 = 1.0$e-$10$ is imposed on the inner solvers.

\begin{figure*}
    \centering
    \begin{subfigure}[b]{0.39\textwidth}
				\includegraphics[width=\textwidth]{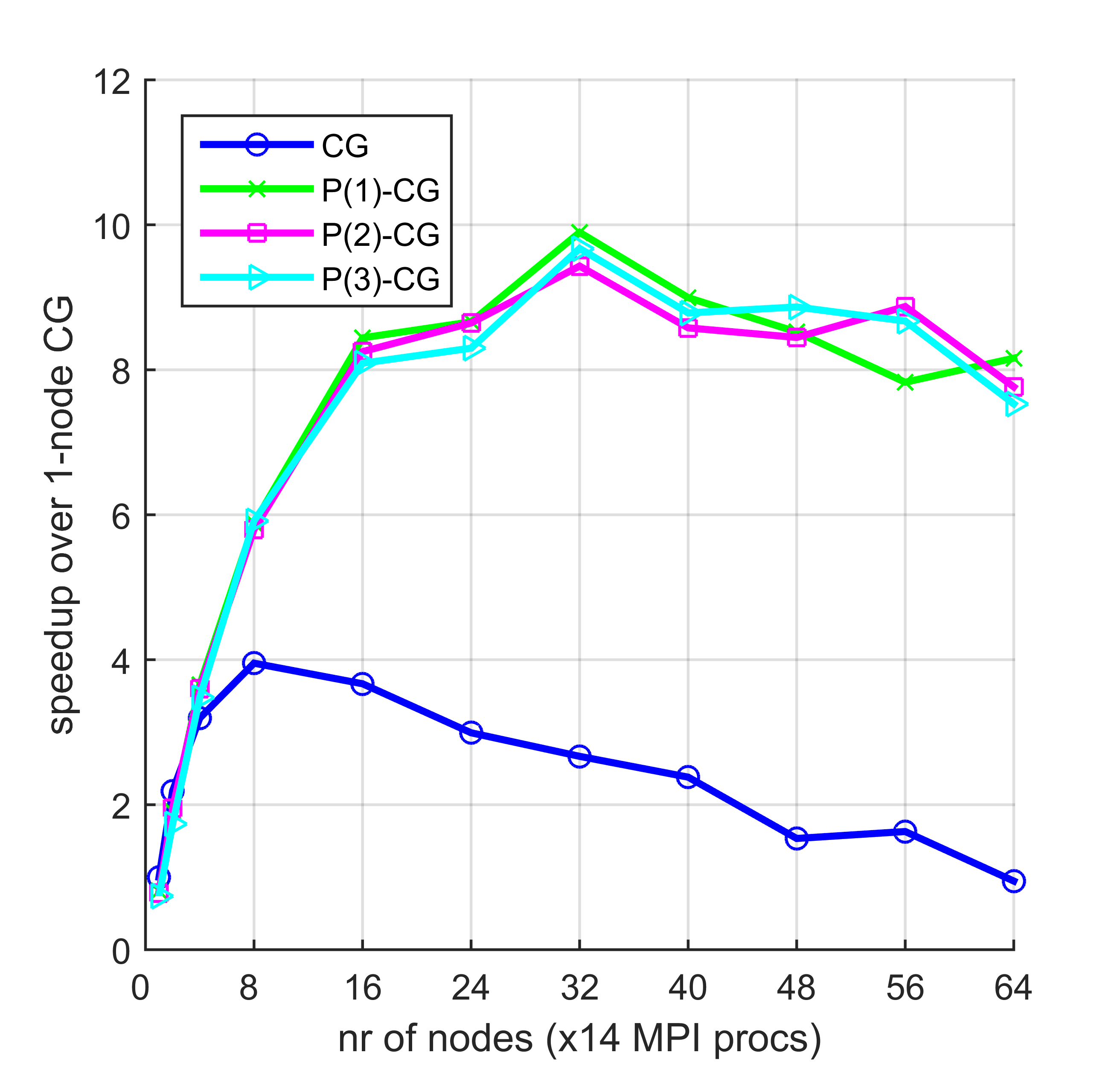} 
        \caption{}
        \label{fig:figure5a}
    \end{subfigure}
		\quad 
    \begin{subfigure}[b]{0.58\textwidth}
				\includegraphics[width=\textwidth]{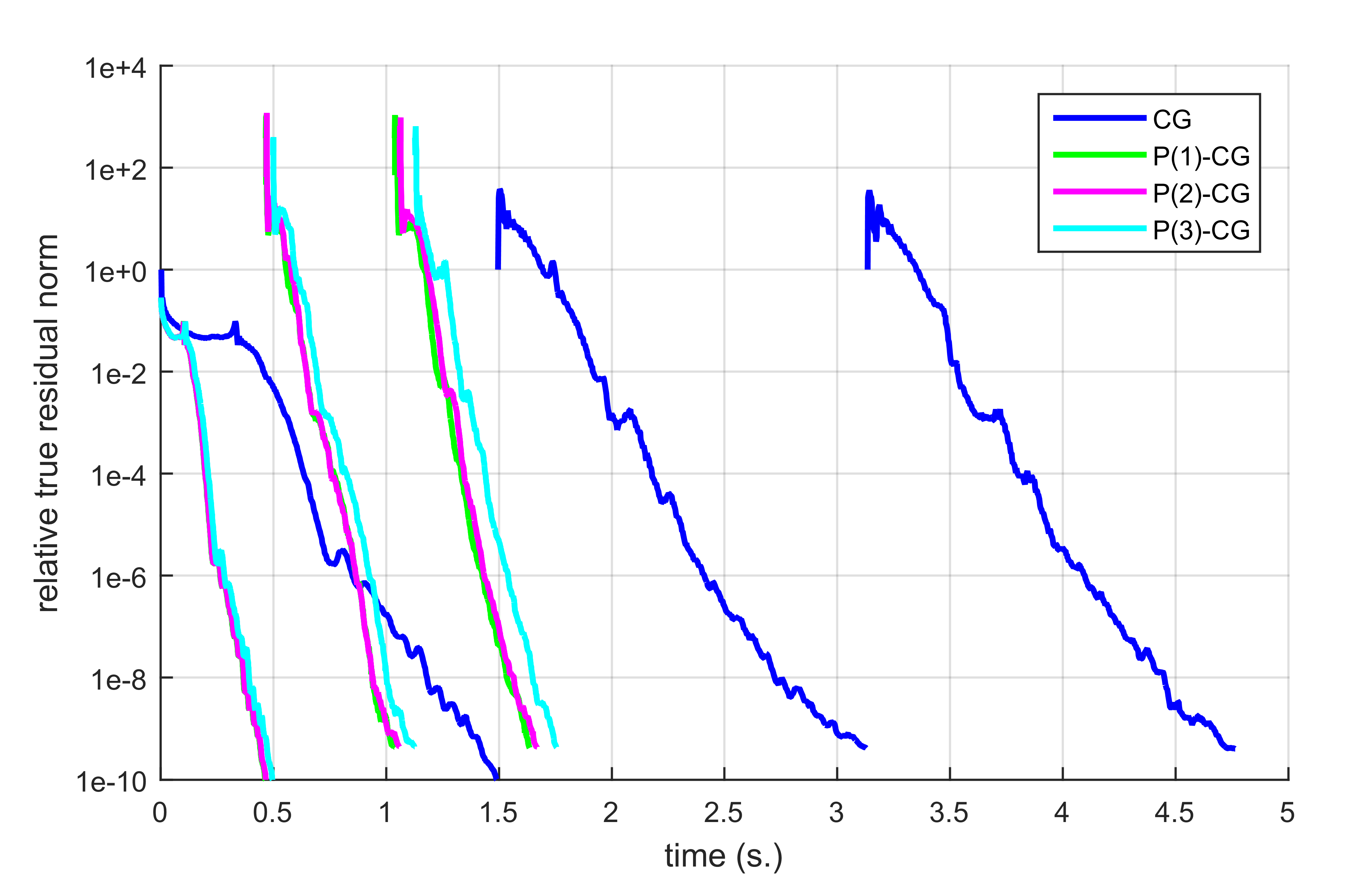} 
        \caption{}
        \label{fig:figure5b}
    \end{subfigure}
		\vspace{-0.2cm}
\caption{\textbf{(B3)} \small
(a) Bratu Solid Fuel Ignition problem (PETSc SNES ex5) with parameter $\lambda_{B} = 6.0$,
$1,000\times1,000$ grid points. SOR-preconditioned Chebyshev preconditioner. 
Inner 
tolerance: 
$\|b-A\bar{x}_i\|_2/\|b\|_2 \leq 1.0$e-$10$.
(b) Accuracy experiment on $24$ nodes ($336$ processes). Residual norm $\|b-A\bar{x}_i\|_2/\|b\|_2$ 
as a function of total time spent (in s.). 
Outer 
tolerance = $1.0$e-$8$.
}
\label{fig:figure5} 
\end{figure*}

Fig.\,\ref{fig:figure5a} presents the relative speedup of the p($l$)-CG methods for $l\in\{1,2,3\}$ compared to classic CG on one node. Timing data for the p-CG method is not included in the figure since it is unable to achieve the desired accuracy; the maximal accuracy attainable by p-CG is $\|b-A\bar{x}_i\|_2/\|b\|_2 = 1.3$e-$9$/$4.2$e-$8$/$6.4$e-$8$ for outer SNES iteration $1/2/3$ respectively. On $32$ nodes the speedup achieved by the pipelined methods over classic CG is approximately $4 \times$. For this small sized problem the p($l$)-CG method does not scale beyond 32 nodes since CPU times are dominated by the $\mathcal{O}(\log_2($\#nodes$))$ behavior of the global reduction phases in this regime. Longer pipelines do not yield a higher speedup compared to pipeline length $1$ due to the cost of the preconditioner, which requires two 
\textsc{spmv}s each iteration. 

Fig.\,\ref{fig:figure5b} shows the relative residual norm as a function of the total time spent for solving the three 
Newton iterations on a $24$ node setup. The timings and residual norms are taken from the same experiment as Fig.\,\ref{fig:figure5a}. The total time-to-solution is divided by the respective number of iterations for the different CG methods to accurately compare them 
for accurate comparison on a single graph. The p($l$)-CG methods 
outperform the classic CG method despite a small increase in 
total number of iterations 
for longer pipeline lengths. The stable p($l$)-CG algorithm reaches the 
relative tolerance 
$1.0$e-$10$ for any choice of the 
parameter $l$. 

\section{Conclusions} \label{sec:conclusions}

This work presents a redesigned algorithmic variant of the $l$-length pipelined Conjugate Gradient method, p($l$)-CG for short. The main improvement over former pipelined CG variants is the 
significantly improved maximal attainable accuracy that is attained by the new algorithm. 
More specifically, it is shown analytically and verified experimentally that the stable p($l$)-CG method attains the same precision on the solution that is attainable by the classic CG method. By introducing intermediate auxiliary bases the propagation of local rounding errors in the recurrence relations is eliminated, allowing for high-precision solution independently of the choice of the pipeline length $l$. The new p($l$)-CG algorithm is elegant in the sense that it has no additional computational overhead and only minor additional storage requirements compared to previous versions of the p($l$)-CG algorithm. The increased stability thus comes without the cost of increased complexity. The improved algorithm effectively replaces former 
(less stable) pipelined CG variants. As such, this work resolves one of the major numerical issues that has restricted the practical usability of pipelined Krylov subspace methods since their initial development. 

Generalizing the stabilization technique proposed in this work to other pipelined methods is 
non-trivial. 
A similar methodology could 
be applied to the existing pipelined GMRES method, p($l$)-GMRES \cite{ghysels2013hiding}. However, practical restrictions limit the viability of this approach. Notably, the full storage of an additional $(l-1)$ auxiliary bases would be required, which can be assumed to be unfeasible for large-scale applications. An $l$-length variant of the pipelined BiCGStab method \cite{cools2017communication} is currently not available, but the proposed technique could be promising for direct application in the context of Bi-Conjugate Gradient methods. 

The research presented in this manuscript provides 
vital advancements towards establishing a numerically stable communication hiding variant of the Conjugate Gradient method. However, it is well-known that Krylov subspace methods may also suffer from delay of convergence due to loss of basis orthogonality. Our experiments indicate that this effect is typically enlarged 
as a function of pipeline length. 
Analyzing 
the stability issues related to loss of orthogonality deserves to be treated as part of 
future work.

\section*{Acknowledgments}
The authors are grateful for the funding received that supported this work. In particular, 
J.\,C.~acknowledges funding by the University of Antwerp Research Council under the University Research Fund (BOF) and
S.\,C.~is funded by the Flemish Research Foundation (FWO Flanders) under grant 12H4617N. 
We also cordially thank Pieter Ghysels for providing useful comments on an earlier version of this manuscript.

\bibliographystyle{plain}
\bibliography{refs2}

\vspace{-0.9cm}

\end{document}